
\documentstyle{amsppt}
\baselineskip18pt
\magnification=\magstep1
\pagewidth{30pc}
\pageheight{45pc}

\hyphenation{co-deter-min-ant co-deter-min-ants pa-ra-met-rised
pre-print pro-pa-gat-ing pro-pa-gate
fel-low-ship Cox-et-er dis-trib-ut-ive}
\def\leaderfill{\leaders\hbox to 1em{\hss.\hss}\hfill}
\def\A{{\Cal A}}

\def\H{{\Cal H}}

\def\L{{\Cal L}}
\def\latl#1{{\Cal L}_L^{#1}}
\def\latr#1{{\Cal L}_R^{#1}}
\def\lat#1{{\Cal L}^{#1}}

\

\def\ldescent#1{{\Cal L (#1)}}
\def\rdescent#1{{\Cal R (#1)}}

\def\idest{i.e.,\ }

\def\d{{\delta}}

\def\e{{\varepsilon}}

\def\th{{\theta}}

\def\l{{\lambda}}

\def\t{{\tau}}

\def\T{{\widetilde T}}
\def\te{\widetilde t}

\def\tmu{\tilde{\mu}}
\def\tem{\tilde{M}}

\def\br{{\bold r}}

\def\b0{\text{\bf 0}}

\def\ra{{\ \longrightarrow \ }}

\def\lan{{\langle}}
\def\ran{{\rangle}}

\def\lrt#1#2{\left\langle {#1}, {#2} \right\rangle}
\def\lrh#1#2{\left\langle {#1}, {#2} \right\rangle_\H}

\def\zed{{\Bbb Z}}
\def\kyu{{\Bbb Q}}

\def\boxit#1{\vbox{\hrule\hbox{\vrule \kern3pt
\vbox{\kern3pt\hbox{#1}\kern3pt}\kern3pt\vrule}\hrule}}
\def\rabbit{\vbox{\hbox{\kern0pt
\vbox{\kern0pt{\hbox{---}}\kern3.5pt}}}}

\def\tableau#1{
        \hbox {
                \hskip -10pt plus0pt minus0pt
                \raise\baselineskip\hbox{
                \offinterlineskip
                \hbox{#1}}
                \hskip0.25em
        }
}

\def\tabCol#1{
\hbox{\vtop{\hrule
\halign{\strut\vrule\hskip0.5em##\hskip0.5em\hfill\vrule\cr\lower0pt
\hbox\bgroup$#1$\egroup \cr}
\hrule
} } \hskip -10.5pt plus0pt minus0pt}

\def\CR{
        $\egroup\cr
        \noalign{\hrule}
        \lower0pt\hbox\bgroup$
}



\def\blank#1#2{
\hbox to #1{\hfill \vbox to #2{\vfill}}
}


\def\strut{\vrule height10pt depth5pt width0pt}


\def\ignore#1{{}}

\def\seca{1}
\def\secb{2}
\def\secc{3}
\def\secd{4}
\def\sece{5}
\def\secf{6}
\def\secg{7}
\def\secconc{8}

\topmatter
\title Generalized Jones traces and Kazhdan--Lusztig bases
\endtitle

\author R.M. Green \endauthor
\affil Department of Mathematics \\ University of Colorado \\
Campus Box 395 \\ Boulder, CO  80309-0395 \\ USA \\ {\it  E-mail:}
rmg\@euclid.colorado.edu \\
\newline
\endaffil

\abstract 
We develop some applications of certain algebraic and combinatorial conditions
on the elements of Coxeter groups, such as elementary proofs of the 
positivity of certain structure constants for the associated Kazhdan--Lusztig 
basis.  We also explore some consequences of the existence of a 
Jones-type trace on the Hecke algebra of a Coxeter group,
such as simple procedures for
computing leading terms of certain Kazhdan--Lusztig polynomials.
\endabstract


\endtopmatter

\centerline{\bf To appear in the Journal of Pure and Applied Algebra}

\head Introduction \endhead

In their seminal paper, Kazhdan and Lusztig \cite{{\bf 18}} defined some remarkable
bases, $\{C_w : w \in W\}$ and $\{C'_w : w \in W\}$ for the Hecke algebra 
$\H$ of
an arbitrary Coxeter group $W$.  The construction of these Kazhdan--Lusztig
bases from the obvious basis $\{T_w : w \in W\}$ of the Hecke algebra 
involves certain polynomials, $\{P_{y, w}(q) : y, w \in W\}$, now known as
Kazhdan--Lusztig polynomials.  When $y < w$ in the Bruhat order on $W$,
$P_{y, w}(q)$ is of degree at most $(\ell(w) - \ell(y) - 1)/2$, where $\ell$
is the length function on the Coxeter group.  The cases where this degree
bound is achieved are of particular importance, and in such cases, the
leading coefficient of $P_{y, w}(q)$ is denoted by $\mu(y, w)$.  The 
$P_{y, w}(q)$ and $\mu(y, w)$ are defined by recurrence relations and are
very difficult to compute efficiently, even for some moderately small groups.

The Kazhdan--Lusztig bases have some remarkable and subtle properties.  One of
these is that (at least in the well-understood cases) if we write $$
C'_x C'_y = \sum_{z \in W} f_{x, y, z} C'_z
,$$ the structure constants $f_{x, y, z}$ are Laurent polynomials with 
nonnegative integer coefficients.  No elementary proof of this phenomenon
has ever been found, except in easy cases such as the dihedral groups
(type $I_2(m)$).  It is, however, possible to establish partial results
in this direction using elementary (i.e., algebraic or combinatorial) means.  
For example, recent work of Geck \cite{{\bf 5}, Theorem 5.10} proves the weaker 
result that positivity of structure constants holds for the asymptotic
Hecke algebra associated to the symmetric group (i.e., Coxeter type $A$).

Like Geck's paper \cite{{\bf 5}}, this paper is motivated by a desire to understand
the Kazhdan--Lusztig bases as far as possible, using elementary methods
and a relatively small set of hypotheses, which themselves should be
verifiable using elementary means.  We aim for conceptual
proofs rather than case by case checks based on Coxeter graphs or the
classification of Kazhdan--Lusztig cells; in particular, we do not restrict
our attention to finite and affine Weyl groups, where the Kazhdan--Lusztig 
theory is best understood.

There are four main hypotheses used in this paper.  The principal one
(Property B)
concerns the existence of a certain remarkable kind of trace on the Hecke
algebra, which we conjecture exists in general.  In type $A$, such a trace
arises by an appropriate scaling of Jones' well-known trace on the
Hecke algebra, and its quotient the Temperley--Lieb algebra
\cite{{\bf 17}, \S11}.  As we will explain, such traces are also known to exist in 
other cases, and they may often be constructed to have the Markov property.
It is possible, although not very easy, to prove the existence of such 
traces in certain special cases by using elementary arguments.  Two of the
other hypotheses that we use (Property F and Property S) are combinatorial 
criteria that are fairly easy to check in particular cases.  The fourth
criterion, Property W, is a weaker, algebraic version of Property S.
All the proofs in the present paper are elementary and are largely 
self-contained.

Our main tool in the present paper is the Kazhdan--Lusztig type basis 
$\{c_w\}$ of the Temperley--Lieb quotient $TL(X)$.  This basis, which is
indexed by the fully commutative 
elements of the Coxeter group, in the sense of Stembridge \cite{{\bf 25}},
was introduced for arbitrary Coxeter groups $W(X)$ by J. Losonczy and
the author in \cite{{\bf 13}}.  

Theorem \sece.13 shows how, in the presence of Property F and Property W, the 
structure constants with respect to the $c$-basis are closely related to 
leading terms of Kazhdan--Lusztig polynomials.

Theorem \secf.13 shows that, under the same hypotheses, the structure constants
for the $c$-basis are Laurent polynomials with nonnegative coefficients.
Theorem \secf.16 shows that if one additionally assumes Property S, it can be
shown that if $z$ is fully commutative, the coefficient of the Kazhdan--Lusztig
basis element $C'_z$ in any product $C'_x C'_y$ is also a nonnegative
Laurent polynomial.

Theorem 7.10 shows how, in the presence of Property B and Property F and
a bipartite Coxeter graph, the leading coefficients $\mu(x, y)$ (where
$x, y$ are fully commutative) can be
computed very easily using suitable traces, assuming these can be explicitly
constructed, which they often can.  This appears to be new even in type $A$,
in which case one can compute the coefficients using Jones' trace from
\cite{{\bf 17}} (see Example 7.15).

Apart from the applications to Kazhdan--Lusztig polynomials and bases, 
our results can be used to bring various 
theorems in the literature into a single context.  As we shall mention, 
many of the results
in the literature on the elements $C'_w$ in the case where $w$ is fully 
commutative are either closely related to the existence of the traces
mentioned above, or are proving additional properties about them in the
cases where they do exist.  The traces are thus of central importance in
the study of these questions.

\head \seca.  Hecke algebras \endhead

Let $X$ be a Coxeter graph, of arbitrary type,
and let $W(X)$ be the associated Coxeter group with distinguished
(finite) set of generating involutions $S(X)$.  In other words, $W = W(X)$ is 
given by the presentation $$
W = \lan S(X) \ | \ (st)^{m(s, t)} = 1 \text{ for } m(s, t) < \infty \ran
,$$ where $m(s, s) = 1$.  (It turns out that the elements 
of $S = S(X)$ are distinct as group elements, and that $m(s, t)$ is the order
of $st$.)  Denote by $\H_q = \H_q(X)$ the Hecke
algebra associated to $W$.  This is a $\zed[q, q^{-1}]$-algebra
with a basis consisting of (invertible) elements $T_w$, with $w$ ranging over 
$W$, satisfying $$T_s T_w = 
\cases
T_{sw} & \text{ if } \ell(sw) > \ell(w),\cr
q T_{sw} + (q-1) T_w & \text{ if } \ell(sw) < \ell(w),\cr
\endcases$$ where $\ell$ is the length function on the Coxeter group
$W$, $w \in W$, and $s \in S$.

For many applications it is convenient to extend the scalars of $\H_q$ to
produce an $\A$-algebra $\H$,
where $\A = \zed[v, v^{-1}]$ and $v^2 = q$, and to define a scaled version
of the $T$-basis, $\{\T_w : w \in W\}$, where $\T_w := v^{-\ell(w)} T_w$.
We will write $\A^+$ and $\A^-$ for $\zed[v]$ and $\zed[v^{-1}]$, respectively,
and we denote the $\zed$-linear ring
homomorphism $\A \ra \A$ exchanging $v$ and $v^{-1}$ by $\bar{\ }$.  We
can extend $\bar{\ }$ to a ring automorphism of $\H$
(as in \cite{{\bf 6}, Theorem 11.1.10}) by the condition that $$
\overline{\sum_{w \in W} a_w \T_w} := \sum_{w \in W} \overline{a_w}
\T_{w^{-1}}^{-1}
,$$ where the $a_w$ are elements of $\A$.

In \cite{{\bf 18}}, Kazhdan and Lusztig proved the following

\proclaim{Theorem \seca.1. (Kazhdan, Lusztig)}
For each $w \in W$, there exists a unique $C'_w \in \H$ such that both
$\overline{C'_w} = C'_w$ and $$
C'_w = \T_w + \sum_{y < w} a_y \T_y
,$$ where $<$ is the Bruhat order on $W$ and $a_y \in v^{-1} \A^-$.
The set $\{C'_w : w \in W\}$ forms an $\A$-basis for $\H$.
\qed\endproclaim

Following \cite{{\bf 6}, \S11.1}, we denote the coefficient of $\T_y$ in $C'_w$ 
by $P^*_{y, w}$.  The {\it Kazhdan--Lusztig polynomial} $P_{y, w}$ is then 
given by $v^{\ell(w) - \ell(y)} P^*_{y, w}$.

\proclaim{Proposition \seca.2}
Define a symmetric $\A$-bilinear 
form, $\lrh{\ }{\ }$, on $\H$ by $$
\lrh{T_x}{T_y} = \d_{x, y} q^{\ell(x)}
,$$ where $\d$ is the Kronecker delta.  Let $x, y \in W$ and $s \in S$. 
\item{\rm (i)}{We have $
\lrh{T_s T_x}{T_y} = \lrh{T_x}{T_s T_y}
,$ and thus $
\lrh{T_x}{T_y} = \lrh{T_x T_{y^{-1}}}{1}
.$  If $*$ denotes the $\zed[q, q^{-1}]$-linear map from $\H_q$ to $\H_q$ 
sending $T_w$ to $T_w^{-1}$, then $\lrh{hh_1}{h_2} = \lrh{h_1}{h^* h_2}$ 
for all $h, h_1, h_2 \in \H_q$.}
\item{\rm (ii)}
{The form $\lrh{\ }{\ }$ induces a nondegenerate trace $
\t_\H : \H \ra \A
$ given by $\t_\H(a) = \lrh{a}{1}$, and we have
$\t_\H(ab) = \t_\H(ba)$ for all $a, b \in \H$.
The restriction of $\t_\H$ to $\H_q$ takes values in $\zed[q, q^{-1}]$.}
\item{\rm (iii)}
{The basis $\{ \T_w : w \in W \}$ is orthonormal with respect to
$\lrh{\ }{\ }$.}
\endproclaim

\demo{Proof}
This is a routine exercise using the definition of $\H$; see
\cite{{\bf 6}, Theorem 8.1.1} for more details.
\qed\enddemo

The following well-known result shows how the form $\lrh{\ }{\ }$ is 
well-suited to studying questions about the Kazhdan--Lusztig basis.

\proclaim{Proposition \seca.3}
\item{\rm (i)}
{The basis $\{C'_w : w \in W\}$ is almost orthonormal
with respect to the form $\lrh{\ }{\ }$: in other words, whenever 
$w, w' \in W$, we have $$
\lrh{C'_w}{C'_{w'}} = 
\cases 1 \mod v^{-1} \A^- & \text{ if } w = w',\cr
0 \mod v^{-1} \A^- & \text{ otherwise.} \cr
\endcases$$}\item{\rm (ii)}
{Suppose $x \in \H$ satisfies both $\bar{x} = x$ and $\lrh{x}{x} = 
1 \mod v^{-1} \A^-$.  Then either $x$ or $-x$ is one of the Kazhdan--Lusztig
basis elements $C'_w$ for some $w$.}
\endproclaim

\demo{Proof}
Part (i) follows easily from Proposition \seca.2 (iii) and Theorem \seca.1.

Part (ii) is a well-known result of Lusztig (compare with \cite{{\bf 21}, Theorem
14.2.3}), which can be proved using similar methods.
\qed\enddemo

\head \secb. Property B and homogeneous traces \endhead

Let $J(X)$ be the two-sided ideal of $\H$ generated by the elements $$
\sum_{w \in \lan s, s' \ran} T_w,
$$ where $(s, s')$ runs over all pairs of elements of $S$
that correspond to adjacent nodes in the Coxeter graph, and $\lan s, s' \ran$
is the parabolic subgroup generated by $s$ and $s'$.  
(If the nodes corresponding to $(s, s')$ are connected by a
bond of infinite strength, then we omit the corresponding relation.)

Following Graham \cite{{\bf 7}, Definition 6.1}, we define the {\it generalized
Temperley--Lieb algebra} $TL(X)$ to be
the quotient $\A$-algebra $\H(X)/J(X)$.  We denote the corresponding
epimorphism of algebras by $\th : \H(X) \ra TL(X)$.  Since the generators
of $J(X)$ lie in $\H_q(X)$, we also obtain a $\zed[q, q^{-1}]$-form
$TL_q(X)$, of $TL(X)$.
Let $t_w$ (respectively, $\te_w$) denote the image in 
$TL(X)$ of the basis element $T_w$ (respectively, $\T_w$) of $\H$.

A product $w_1w_2\cdots w_n$ of elements $w_i\in W$ is called
{\it reduced} if \newline $\ell(w_1w_2\cdots w_n)=\sum_i\ell(w_i)$.  We reserve
the terminology {\it reduced expression} for reduced products 
$w_1w_2\cdots w_n$ in which every $w_i \in S$.  We write $$
\ldescent{w} = \{s \in S : \ell(sw) < \ell(w)\}
$$ and $$
\rdescent{w} = \{s \in S : \ell(ws) < \ell(w)\}
.$$  The set $\ldescent{w}$ (respectively, $\rdescent{w}$) is called the 
{\it left} (respectively, {\it right}) {\it descent set} of $w$.

Call an element $w \in W$ {\it complex} if it can be written 
as a reduced product $x_1 w_{ss'} x_2$, where $x_1, x_2 \in W$ and
$w_{ss'}$ is the longest element of some rank 2 parabolic subgroup 
$\lan s, s'\ran$ such that $s$ and $s'$ correspond to adjacent nodes
in the Coxeter graph.
Denote by $W_c(X)$ the set of all elements of $W$
that are not complex.  The elements of $W_c = W_c(X)$ are the 
{\it fully commutative}
elements of \cite{{\bf 25}}; they are characterized by the property that any two
of their reduced expressions may be obtained from each other by repeated
commutation of adjacent generators.

We define the $\A^-$-submodule $\L$ of
$TL(X)$ to be that generated by the $\{\te_w : w \in W_c\}$.  We define
$\pi : \L \ra \L/v^{-1}\L$ to be the canonical $\zed$-linear projection.

By \cite{{\bf 13}, Lemma 1.4}, the ideal $J(X)$ is fixed by $\bar{\ }$, so 
$\bar{\ }$ induces an involution on $TL(X)$, which we also denote by 
$\bar{\ }$.

The next result is an analogue of Theorem \seca.1, and the proof is 
similar; in particular, it works for arbitrary Coxeter groups.  The basis
elements $\{c_w : w \in W_c\}$ may be regarded as baby versions of the 
Kazhdan--Lusztig basis elements $C'_w$.  We will prove in Proposition 
\secf.3 (i) below that, under certain hypotheses, we have $\theta(C'_w)
= c_w$ if $w \in W_c$.  This is the eponymous ``projection property'' of
\cite{{\bf 14}}.  There is no known example of a Coxeter group that fails
to satisfy this projection property.  Although it
is not generally true that $\theta(C'_w) = 0$ for $w \not\in W_c$,
many Coxeter groups do have this latter property, such as those of type 
$A_n$, $B_n$,
$F_4$, $H_3$, $H_4$, $I_2(m)$, $\widehat{A}_n$ and $\widehat{C}_n$.  We
will discuss this in detail later; see, for example, the remarks following
Theorem \secf.13.

\proclaim{Theorem \secb.1}
\item{\rm (i)}
{The set $\{t_w : w \in W_c\}$ is a $\zed[q, q^{-1}]$-basis for $TL_q(X)$.
The set $\{\te_w : w \in W_c\}$ is an $\A$-basis for $TL(X)$, and an
$\A^-$-basis for $\L$.}
\item{\rm (ii)}
{For each $w \in W_c$, there exists a unique $c_w \in TL(X)$ such that both
$\overline{c_w} = c_w$ and $\pi(c_w) = \pi(\te_w)$.  Furthermore, we have $$
c_w = \te_w + \sum_{{y < w} \atop {y \in W_c}} a_y \te_y
,$$ where $<$ is the Bruhat order on $W$.}
\item{\rm (iii)}
{The set $\{c_w : w \in W_c\}$ forms an $\A$-basis for $TL(X)$ and an
$\A^-$-basis for $\L$.}
\item{\rm (iv)}
{If $x \in \L$ and $\bar{x} = x$, then $x$ is a
$\zed$-linear combination of the $c_w$.}
\item{\rm (v)}
{There is an $\A$-linear anti-automorphism, $*$, of $TL(X)$ that sends
$\te_w$ to $\te_{w^{-1}}$ and $c_w$ to $c_{w^{-1}}$ for all $w \in W_c$.}
\endproclaim

\demo{Proof}
Part (i) is \cite{{\bf 7}, Theorem 6.2}, and 
parts (ii) and (iii) are \cite{{\bf 13}, Theorem 2.3}, except for the 
assertions about $\L$, which are
are immediate from the definitions.  Part (iv) follows from
(ii) and the fact that $$
\overline{\sum_{u \in W_c} a_u c_u} = \sum_{u \in W_c} \overline{a_u} c_u
.$$

For part (v), we note that it is well known that the $\zed[q, q^{-1}]$-linear
map from $\H_q$ to $\H_q$ that sends $T_w$ to $T_{w^{-1}}$ is an 
anti-automorphism, $*$, of $\H_q$.  By extending scalars, we obtain an 
$\A$-linear anti-automorphism (also denoted by $*$) of $\H$ that 
sends $\T_w$ to $\T_{w^{-1}}$; furthermore, $*$ commutes with the ring
automorphism $\bar{\ }$.  It is clear
from the definition of $J(X)$ that $J(X)$ is fixed by this map, so we obtain
an anti-automorphism of $TL(X)$ sending $\te_w$ to $\te_{w^{-1}}$, in 
particular, when $w \in W_c$.  Since $*$ and $\bar{\ }$ commute, part (ii)
shows that $*$ sends $c_w$ to $c_{w^{-1}}$.
\qed\enddemo

The following hypothesis is analogous to Proposition \seca.2.

\proclaim{Hypothesis \secb.2}
Let $X$ be an arbitrary Coxeter graph.
There exists a symmetric $\A$-bilinear 
form, $\lrt{\ }{\ }$, on $TL(X)$ satisfying the following properties for
all $x, y \in W_c$ and $s \in S$:
\item{\rm (i)}{$\lrt{\te_s \te_x}{\te_y} = \lrt{\te_x}{\te_s \te_y}$
(and therefore $\lrt{h h_1}{h_2} = \lrt{h_1}{h^* h_2}$ for all $h, h_1,
h_2 \in TL(X)$);} 
\item{\rm (ii)}{the basis $\{\te_w : w \in W_c\}$ is almost 
orthonormal with respect to $\lrt{\ }{\ }$, meaning that $$
\lrt{\te_x}{\te_y} = 
\cases 1 \mod v^{-1} \A^- & \text{ if } x = y,\cr
0 \mod v^{-1} \A^- & \text{ otherwise.} \cr
\endcases$$}
\endproclaim

An immediate consequence of Hypothesis \secb.2 (ii) is that the bilinear form 
$\lrt{\ }{\ }$ restricts to an $\A^-$-valued $\A^-$-form on $\L$.

\definition{Definition \secb.3 (Property B)}
If Hypothesis \secb.2 holds for the Coxeter graph $X$, we say that $X$ (or
$W(X)$) has {\it Property B}.
\enddefinition

Some immediate consequences of Property B are the following.

\proclaim{Proposition \secb.4}
Assume that $W$ has Property B.
\item{\rm (i)}
{The basis $\{c_w : w \in W_c\}$ is almost orthonormal
with respect to the form $\lrt{\ }{\ }$: in other words, whenever 
$x, y \in W_c$, we have $$
\lrt{c_x}{c_{y}} = 
\cases 1 \mod v^{-1} \A^- & \text{ if } x = y,\cr
0 \mod v^{-1} \A^- & \text{ otherwise.} \cr
\endcases$$}\item{\rm (ii)}
{Suppose $x \in TL(X)$ satisfies both $\bar{x} = x$ and $\lrt{x}{x} = 
1 \mod v^{-1} \A^-$.  Then either $x$ or $-x$ is one of the canonical
basis elements $c_w$ for some $w \in W_c$.}
\item{\rm (iii)}
{The form $\lrt{\ }{\ }$ induces a nondegenerate trace 
$\t : TL(X) \ra \A$ given by $\t(a) = \lrt{a}{1}$.  We have 
$\t(ab) = \t(ba)$ for all $a, b \in TL(X)$ and $\t(a^*) = \t(a)$.}
\endproclaim

\demo{Proof}
Part (i) is immediate from Theorem \secb.1 (ii) and Hypothesis \secb.2 (ii).
(In fact, this shows that Hypothesis \secb.2 (ii) and Proposition \secb.4 (i) 
are equivalent.)

Part (ii) is proved by a standard argument, given in 
\cite{{\bf 10}, Proposition 4.3.4}.

We now turn to part (iii).  Symmetry of the form $\lrt{\ }{\ }$ shows that
$\t(ab) = \t(ba)$ for all $a, b \in TL(X)$.
Repeated applications of Hypothesis \secb.2 (i) 
show that $\lrt{\te_w}{1} = \lrt{1}{\te_{w^{-1}}}$, and symmetry of
$\lrt{\ }{\ }$ together with $\A$-bilinearity then show that 
$\t(a^*) = \t(a)$.  Hypothesis \secb.2 (ii) shows that $\lrt{\ }{\ }$ is
nondegenerate, from which it is clear that the associated trace is
nondegenerate.
\qed\enddemo

The main focus of this paper is to explore further consequences of Property B.
Hypothesis 2.2 may be checked combinatorially in special cases, although this
is not easy and one needs to know a lot about the structure of the algebra
$TL(X)$ in order to do this.  Conversely, in the cases where Hypothesis 2.2
is known to hold, we will see later, in the main results, that one can 
deduce information about 
Kazhdan--Lusztig polynomials and structure constants that would otherwise 
be hard to prove.  This stands in contrast to the analogous situation
concerning $\H$ and $\lrh{\ }{\ }$, where questions involving the
Kazhdan--Lusztig basis often turn out to be combinatorially very difficult 
or intractable.

\remark{Remark \secb.5}
Property B is known to be true in various special cases, including
the following.
\item{(i)}
{For Coxeter systems of type $A_n (n \geq 1)$ , $B_n (n \geq 2)$, the extended
$H_n$ series (for arbitrary $n \geq 4$) and the dihedral case of $I_2(m)$,
the hypothesis was proved to hold in \cite{{\bf 10}, Corollary 4.3.3} using
constructive methods from diagram algebras and planar algebras.  In 
particular, one can construct a bilinear
form in type $A$ by $\lrt{c_x}{c_y} = \t(c_x c_{y^{-1}})$, where
$\t$ is obtained from the {\it Jones trace} \cite{{\bf 17}, (11.5)} after 
multiplication by the factor $v^{-(n+1)} (v+v^{-1})^{n+1}$; see also
\cite{{\bf 10}, Definition 3.2.1}.
Note that the symbol $\t$ in \cite{{\bf 17}} corresponds to $(v + v^{-1})^{-2}$
in our notation, and $t$ in \cite{{\bf 17}} corresponds to $v^2$.}
\item{(ii)}
{For Coxeter systems of type $D_n (n \geq 4)$ and the extended $E_n$
series (for arbitrary $n \geq 6$), the hypothesis holds.  Although this is a
consequence of \cite{{\bf 13}, Theorem 3.6} and \cite{{\bf 9}, Theorem 4.3.5},
the proof in \cite{{\bf 9}, \S4.3} that the bilinear form is symmetric contains
a gap.  If the Coxeter group is finite, the argument is completed by
\cite{{\bf 6}, Corollary 8.2.6 (c)}, which shows that any trace $\phi$ on 
$\H$ takes equal values on $T_w$ and $T_{w^{-1}}$, for any $w \in W$.
This gap is also fixable for the cases $E_n$, $n > 8$, or alternatively 
one may describe a trace
satisfying the required property by requiring that whenever $w$ is a reduced
product of $a$ commuting Coxeter generators, we have $$
\t(c_w) = v^{-n} (v+v^{-1})^{n - a}
.$$  As in type $A$, this may be proved using calculi of diagrams: 
the paper \cite{{\bf 8}} describes a diagram calculus for $TL(D_n)$ and 
\cite{{\bf 3}} describes
a (more complicated) diagram calculus for $TL(E_n)$.  Full details of
these constructions will appear in \cite{{\bf 12}}.}
\endremark

We conjecture that Property B holds for all Coxeter groups.

\remark{Remark \secb.6}
Of course, Property B may be reformulated as a conjecture about a
degenerate bilinear form on $\H$ whose radical is precisely $J(X)$.
\endremark

For many of our later purposes, we wish to work with traces $\t$ that are
compatible with the $\zed[q, q^{-1}]$-form of the algebras.

\definition{Definition \secb.7}
Let $A_q$ be an $\zed[q, q^{-1}]$-algebra, and let $A = \A 
\otimes_{\zed[q, q^{-1}]} A_q$.  Let $\t : A \ra \A$ be an $\A$-linear map.
We say that $\t$ is {\it homogeneous} if the restriction, $\t_q$, of $\t$
to $A_q$ takes values in $\zed[q, q^{-1}]$.
\enddefinition

\proclaim{Lemma \secb.8}
Suppose that the Coxeter graph $X$ has Property B, and let $\t$
be the trace corresponding to the bilinear form $\lrt{\ }{\ }$.  Then
Hypothesis \secb.2 is also satisfied by a bilinear form whose trace is
homogeneous.
\endproclaim

\demo{Proof}
Let $\t$ be any trace satisfying Hypothesis \secb.2, and let $\t_q$ be its
restriction to $TL_q(X)$.  (Note that $\t_q$ need not take values in 
$\zed[q, q^{-1}]$.)  Let $p : \A \ra \zed[q, q^{-1}]$ be the 
$\zed$-linear map such that $$
p(v^n) = \cases
v^n & \text{ if } n \text{ is even},\cr
0 & \text{ if } n \text{ is odd}.\cr
\endcases$$  Since $TL_q(X)$ is a $\zed[q, q^{-1}]$-algebra, it follows that 
$p \circ \t_q$ is a $\zed[q, q^{-1}]$-valued trace on $TL_q$.  By extending
scalars to $\A$, $p \circ \t_q$ induces a homogeneous trace, $\t'$, on
$TL(X)$, and it is not hard to check that it has the required properties.
\qed\enddemo

\definition{Definition \secb.9}
We call a trace for $TL(X)$ (or its inflation to $\H$) a {\it homogeneous
trace} (or {\it generalized
Jones trace}) if both (a) it corresponds to a bilinear form satisfying
Hypothesis \secb.2 and (b) it is homogeneous in the sense of Definition
\secb.7.  If the form $\lrt{\ }{\ }$ appearing in Hypothesis \secb.2 is
associated to a homogeneous trace, we call $\lrt{\ }{\ }$ a {\it homogeneous
bilinear form}.  If in addition, a homogeneous trace $\t$ satisfies $\t(c_w)
\in {\Bbb Z}^{\geq 0}[v, v^{-1}]$ for all $w \in W_c$, we say that the trace is
{\it positive}.
\enddefinition

All the traces described in Remark \secb.5 may be easily checked to be 
homogeneous and positive.

\head \secc.  Star reducibility, Property F and Property S \endhead

A key concept for this paper is that of a star operation.  These were
introduced in the simply laced case in \cite{{\bf 18}, \S4.1}, and in general in
\cite{{\bf 20}, \S10.2}.

\definition{Definition \secc.1}
Let $W$ be any Coxeter group and let $I = \{s, t\} \subseteq S$ be a
pair of noncommuting generators whose product has order $m$ (where $m = \infty$
is allowed).  Let $W^I$ denote the set of all $w \in W$ satisfying 
$\ldescent{w} \cap I = \emptyset$.  Standard properties of Coxeter groups 
\cite{{\bf 16}, \S5.12} show that any element $w \in W$ may be uniquely
written as $w = w_I w^I$, where $w_I \in W_I = \lan s, t \ran$
and $\ell(w) = \ell(w_I) + \ell(w^I)$.  
There are four possibilities for elements $w \in W$:
\item{(i)}{$w$ is the shortest element in the coset $W_I w$, so $w_I = 1$ and 
$w \in W^I$;}
\item{(ii)}{$w$ is the longest element in the coset $W_I w$, so $w_I$ is 
the longest element of $W_I$ (which can only happen if $W_I$ is finite);}
\item{(iii)}{$w$ is one of the $(m-1)$ elements $sw^I$, $tsw^I$, $stsw^I, 
\ldots$;}
\item{(iv)}{$w$ is one of the $(m-1)$ elements $tw^I$, $stw^I$, $tstw^I, 
\ldots$.}

The sequences appearing in (iii) and (iv) are called {\it (left) 
$\{s, t\}$-strings},
or {\it strings} if the context is clear.  If $x$ and $y$ are two elements 
of an $\{s, t\}$-string such that $\ell(x) = \ell(y) - 1$, we call the pair 
$\{x, y\}$ {\it left $\{s, t\}$-adjacent}, and we say that $y$ is 
{\it left star reducible to $x$}.

The above concepts all have right-handed counterparts, leading to the notion of
{\it right $\{s, t\}$-adjacent} and {\it right star reducible} pairs of 
elements, and coset decompositions $({^Iw})({_Iw})$.

If there is a (possibly trivial) sequence $$
x = w_0, w_1, \ldots, w_k = y
$$ where, for each $0 \leq i < k$, $w_{i+1}$ is left star reducible or
right star reducible to $w_i$ 
with respect to some pair $\{s_i, t_i\}$, we 
say that $y$ is {\it star reducible to $x$}.  Because star reducibility
decreases length, it is clear that this defines a partial order on $W$.

If $w$ is an element of an $\{s, t\}$-string, $S_w$, we have 
$\{\ell(sw), \ell(tw)\}$ = $\{\ell(w) - 1, 
\ell(w) + 1\}$; let us assume without loss of generality that $sw$ is longer
than $w$ and $tw$ is shorter.  If $sw$ is an element of $S_w$, we define
$^*w = sw$; if not, $^*w$ is undefined.  If $tw$ is an element of $S_w$,
we define $_*w = tw$; if not, $_*w$ is undefined.

There are also obvious right handed analogues to the above concepts,
so the symbols $w^*$ and $w_*$ may be used with the analogous meanings.
\enddefinition

\example{Example \secc.2}
In the Coxeter group of type $B_2$ with $w = ts$, we have $$
{_*w} = s, \ {^*w} = sts, \ w_* = t \text{\ and\ } w^* = tst
.$$  If $x = sts$ then $^*x$ and $x^*$ are undefined; if $x = t$ then
$_*x$ and $x_*$ are undefined.
\endexample

Star reducibility allows us to give concise definitions of the two main
combinatorial criteria of interest in this paper.

\definition{Definition \secc.3 (Property F)}
We say that a Coxeter group $W(X)$, or its Coxeter graph $X$, has
{\it Property F} if every element of $W_c$ is star reducible to a
product of commuting generators from $S$.
\enddefinition

\definition{Definition \secc.4 (Property S)}
We say that a Coxeter group $W(X)$, or its Coxeter graph $X$, has
{\it Property S} if every element of $W(X) \backslash W_c$ is 
star reducible to an element $w$ for which either $\ldescent{w}$ or 
$\rdescent{w}$
(or both) contains a pair of noncommuting generators.
\enddefinition

\remark{Remark \secc.5}
Property F is so called because it is a restatement of the notion of
cancellability which arises in the work of Fan \cite{{\bf 4}}.  
The argument of \cite{{\bf 4}, Lemma 4.3.1} combined with \cite{{\bf 25}, Proposition
2.3} shows that Property F holds
for all Coxeter groups $W$ for which $W_c$ is finite; such groups
were classified independently by Graham \cite{{\bf 7}} and Stembridge \cite{{\bf 25}},
and the connected components of their Coxeter graphs fall into seven infinite 
families: $A$, $B$, $D$, $E$, $F$, $H$ and $I$.  (This is a superset of the 
classification of finite Coxeter groups, but with extended $E_n$, $F_n$ and
$H_n$ series.)  

Property F is not true for arbitrary Coxeter groups, but it does hold
in some other cases.  These include type $\widehat{A}_n$ for $n$ even,
type $\widehat{C}_n$ for $n$ even, type $\widehat{E}_6$, and the case where
$X$ is obtained from the graph of type $A_6$ by relabelling the middle edge
with $4$.  A complete classification for finitely generated Coxeter groups
appears in \cite{{\bf 11}, Theorem 6.3}.
\endremark

\remark{Remark \secc.6}
Property S is so called because it is closely related to a criterion 
appearing in the work of Shi \cite{{\bf 22}, {\bf 23}}.  Shi shows \cite{{\bf 23}, Lemma 2.2}
that this holds for any connected, nonbranching Coxeter graph of a finite or 
affine Weyl group, except type $\widehat{F}_4$.  However, the criterion 
fails for Coxeter systems having a parabolic subsystem of type $D_4$: if
the Coxeter generators are numbered $s_1, \ldots, s_4$ so that $s_2$ fails to
commute with the other three generators, then $$
w = s_1 s_3 s_4 s_2 s_1 s_3 s_4
$$ provides a counterexample to Property S.
\endremark

Unlike Property B, Properties F and S can typically be checked in specific
cases by using fairly short elementary arguments.  As one might guess from
the formulations of these two properties, they complement each other to
some extent and our strongest results are obtained when both properties hold.

The following lemma is extremely useful in inductive arguments.

\proclaim{Lemma \secc.7}
Let $W$ be a Coxeter group with Property B, and let $s, t \in S$ be
noncommuting generators.  Let $x, y \in TL(X)$.  Then we have
\item{\rm (i)}{$
\lrt{\te_s \te_t x}{\te_s y} 
= 
\lrt{\te_t x}{\te_t \te_s y} 
+ \lrt{x}{\te_s y} 
- \lrt{\te_t x}{y} 
;$} \item{\rm (ii)}{$
\lrt{x \te_t \te_s}{y \te_s} 
= 
\lrt{x \te_t}{y \te_s \te_t} 
+ \lrt{x}{y \te_s} 
- \lrt{x \te_t}{y} 
.$}
\endproclaim

\demo{Proof}
Part (i) is immediate from Hypothesis \secb.2 (i) and the identity $$
\T_s \T_s \T_t - \T_t = \T_s \T_t \T_t - \T_s
$$ in $\H$, and part (ii) follows similarly.
\qed\enddemo

\head \secd.  The $\A^-$-lattice $\L$ and Property W \endhead

In this section, we develop some important properties of the $\A^-$-module $\L$
from \S\secb.  The following standard result will be used freely in the sequel.

\proclaim{Lemma \secd.1}
Suppose $W$ has Property B, let $a \in TL(X)$ and let $w \in W_c$.  
If $a \in \L$, the coefficient of $\te_w$ (respectively,
$c_w$) in $a$ with respect to the $\te$-basis (respectively, the $c$-basis)
is equal modulo $v^{-1} \A^-$ to both $\lrt{a}{\te_w}$ and $\lrt{a}{c_w}$.
\endproclaim

\demo{Proof}
This is an immediate consequence of almost orthonormality (see 
Hypothesis \secb.2 (ii) and Proposition \secb.4 (i)).
\qed\enddemo

\definition{Definition \secd.2}
An element $w \in W$ is said to be {\it weakly complex} if (a) it is
complex (in the sense of \S\secb) and (b) it is of the form $w = su$,
where $s \in S$ and $u$ is not complex.
Note that, with the above notation, it must be the case that $su > u$.
\enddefinition

The following definition will be a very useful hypothesis in various results
in the sequel.

\definition{Definition \secd.3 (Property W)}
We say the Coxeter group $W$ has {\it Property W} if, whenever $x \in W$
is weakly complex, we have $\te_x \in v^{-1} \L$.
\enddefinition

\remark{Remark \secd.4}
We shall see in Corollary \secf.15 below  that Property S implies Property W.
In fact, Property F implies Property W (see \cite{{\bf 11}, Theorem
4.6 (i)} for a proof), but this requires much more work than
Proposition \secd.12 below.
Property W seems to be subtle, and is typically difficult to verify or
refute in the absence of any of the aforementioned stronger properties.
We do not know of an example of a Coxeter group that fails to have Property W.
\endremark


\proclaim{Lemma \secd.5}
Let $W$ be any Coxeter group, let $w \in W_c$ and $s \in S$, and 
suppose that $sw \not\in W_c$, in other words, that $sw$ is weakly complex.  
\item{\rm (i)}
{We have $w = w_1 w_2 w_3$ reduced, where (a) every generator
occurring in $w_1$ is distinct from $s$ and commutes with $s$, and (b)
$w_2$ is an alternating product $tsts\ldots$ of length $m(s, t) - 1$, where
$m(s, t)$ is the order of $st$.  It follows that $sw$ has reduced
expressions of the form $s w_1 w_2 w_3$ and $w_1 sw_2 w_3$.}
\item{\rm (ii)}
{If $w \in W_c$ and $u \in S$, then $uw < w \Rightarrow uw \in W_c$, and 
$wu < w \Rightarrow wu \in W_c$.}
\item{\rm (iii)}
{If $u \in S$ and $y \in W$ is such that we have either $w = uy$ or $w = yu$
reduced, then either $sy \in W_c$ or $sy$ is weakly complex.}
\endproclaim

\demo{Proof}
Part (i) is a consequence of \cite{{\bf 25}, Proposition 2.3}, and part (ii) is
immediate from the definition of $W_c$.

For part (iii), note that $y \in W_c$ by (ii).  If $sy < y$ then 
$sy \in W_c$ by (ii).
If $sy > y$ and $sy \not\in W_c$, then $sy$ is weakly complex by definition.
\qed\enddemo

Note that there is an obvious right handed version of Lemma \secd.5.

The following simple result was stated for simply laced Coxeter groups
in \cite{{\bf 7}, Proposition 9.14 (i)} (see also \cite{{\bf 24}, Proposition 2.10}).

\proclaim{Lemma \secd.6}
Let $W$ be an arbitrary Coxeter group and $w \in W_c$.  Let $I = \{s, t\}$
be a pair of noncommuting generators, and take star operations with respect
to $I$.  If $^*w$
(respectively, ${_*w}$, $w^*$, $w_*$) is defined, then $^*w$
(respectively, ${_*w}$, $w^*$, $w_*$) lies in $W_c$.
\endproclaim

\demo{Proof}
The cases of ${_*w}$ and $w_*$ are easy to deal with.  Applying Lemma \secd.5,
we see that if ${^*w} \not \in W_c$, then ${^*w}$ has a reduced expression
beginning with $w_{st}$.  This means
that ${^*w}$ does not lie in the required $\{s, t\}$-string, a contradiction.
The case of $w^*$ follows by a symmetrical argument.
\qed\enddemo

\proclaim{Lemma \secd.7}
Maintain the notation of Lemma \secd.5, and denote by $w'_2$ 
the unique element of $W$ such that $\{w_2, w'_2\}$ are the two elements of
$\lan s, t \ran$ with length $m(s, t) - 1$.  Let us write $$
\te_{sw} = \sum_{u \in W_c} a_u \te_u
.$$  If $a_u \ne 0$, then we have the following:
\item{\rm (i)}{$\ell(u) \leq \ell(w)$;}
\item{\rm (ii)}{we can only have $\ell(u) = \ell(w)$ if 
$u = w_1 w_2 w_3 = w$ or $u = w_1 w'_2 w_3$, and the latter can only occur if 
$w_1 w'_2 w_3$ is an element of $W_c$ of length $\ell(w)$; furthermore,
if $\ell(u) = \ell(w)$, then $a_u = - v^{-1}$.}
\endproclaim

\demo{Proof}
Recall that $TL(X)$ is obtained from $\H$ by the adding the relations $$
\sum_{w \in \lan s, s' \ran} t_w = 0
$$ whenever $\{s, s'\}$ is a pair of noncommuting Coxeter generators generating
a finite (parabolic) subgroup.  Denoting the longest element of this subgroup
by $w_{ss'}$, we can rewrite the relation as $$\eqalignno{
\te_{w_{ss'}} &= - \sum_{w  \in \lan s, s' \ran, w < w_{ss'}}
v^{\ell(w) - \ell(w_{ss'})} \te_w. & (1)
}$$  Using this relation and the other Hecke algebra relations repeatedly,
any element $\te_x (x \in W\backslash W_c)$ can be expressed as a linear 
combination of basis elements $\{\te_u : u \in W_c\}$, where $u < x$.
The assertions now follow from repeated applications of
this relation and Lemma \secd.5 (i).  (The circumstances of (ii) 
can only occur if the relation is applied precisely once.)
\qed\enddemo

For later purposes, it is convenient to define various sublattices of the 
$\A^-$-lattice $\L$.

\definition{Definition \secd.8}
Let $W' \subset W_c$.  We define $\lat{W'}$ to be the free $\A^-$-module with
basis $$
\{ \te_w : w \in W' \} \cup \{ v^{-1} \te_w : w \in W_c \backslash W' \}
.$$  
If $s, t \in S$ are noncommuting generators, 
$W_1 = \{ w \in W_c : sw < w \}$ and $W_2 = \{ w \in W_c : w = stu
\text{ reduced} \}$, we write
$\latl{s}$ and $\latl{st}$ for $\lat{W_1}$ and $\lat{W_2}$, respectively.

One can also define right handed versions, $\latr{s}$ and
$\latr{ts}$, of the above concepts.  Note also that by Theorem \secb.1 (ii),
one can define all these $\A^-$-lattices using the $c$-basis instead of
the $\te$-basis.
\enddefinition

\proclaim{Lemma \secd.9}
Suppose the Coxeter group $W$ has Property B, and let $s \in S$ and
$w \in W_c$ be such that $x = sw$ is weakly complex.  Then we have $\te_x \in 
\latl{s}$.
\endproclaim

\demo{Proof}
Write $x = sw_1 w_2 w_3 = sw$, as in Lemma \secd.5.
The proof is by induction on $\ell(x)$, the case $\ell(x) = 0$ being vacuous.
Let $u \in W_c$.  We need to show that the coefficient of $\te_u$ in
$\te_{sw}$ lies in $\A^-$.  By Lemma \secd.7 (i), we may assume that $\ell(u)
\leq \ell(w)$, \idest that $\ell(u) < \ell(x)$.

If $\ell(u) = \ell(w)$, Lemma \secd.7 (ii) shows that the coefficient of 
$\te_u$ in $\te_{sw}$ is $-v^{-1}$, which satisfies the hypotheses.

Suppose that $\te_x \not\in \L$.  We claim that there exists $y \in W_c$
with $\ell(y) < \ell(w)$ and $\lrt{\te_x}{\te_y} \not\in \A^-$.  By
the assumption, there 
exists an $n > 0$ such that $a = v^{-n} \te_x \in \L$ but $v^{-(n-1)} 
\te_x \not\in \L$, and there exists $y \in W_c$ such that $\te_y$ occurs
with nonzero coefficient in $\te_x$ and such that the coefficient of 
$\te_y$ in $a$ lies in $\A^- \backslash v^{-1} \A^-$.  (By the previous
paragraph, this cannot happen unless $\ell(y) < \ell(w)$.)  Using Lemma
\secd.1, we find that the constant coefficient of $\lrt{a}{\te_y}$ is nonzero,
which means that the coefficient of $v^n$ in $\lrt{\te_x}{\te_y}$
is nonzero, as claimed.  This means that to show that $\te_x \in \L$, it is
sufficient to verify that $\lrt{\te_x}{\te_y} \in \A^-$ when $y \in W_c$
and $\ell(y) < \ell(w)$.  
Assume from now on that $u$ satisfies these properties.

By Property B, we have $$
\lrt{\te_x}{\te_u} = \lrt{\te_s \te_w}{\te_u}
= \lrt{\te_w}{\te_s \te_u}
.$$  There are three subcases to consider.

The first possibility is that $su < u$, in which case we have $$\eqalign{
\lrt{\te_w}{\te_s \te_u} &= \lrt{\te_w}{\te_{su} + (v - v^{-1}) \te_u} \cr
&= \lrt{\te_w}{\te_{su}} + 
(v - v^{-1}) \lrt{\te_w}{\te_u}.\cr
}$$  Since $su \in W_c$ and $\ell(su) < \ell(w)$, Hypothesis \secb.2 (ii) 
shows that 
$\lrt{\te_w}{\te_{su}} \in v^{-1} \A^-$.  Similarly, since $u \in W_c$ and
$\ell(u) < \ell(w)$, we have $\lrt{\te_w}{\te_u} \in v^{-1} \A^-$, and thus
$(v - v^{-1}) \lrt{\te_w}{\te_u} \in \A^-$.

The second possibility is that $su > u$ and $su \in W_c$.  In this case,
we cannot have $su = w$, because $sw > w$ and $s(su) < su$.  Hypothesis
\secb.2 (ii) applies again to show that $\lrt{\te_w}{\te_s \te_u} \in 
v^{-1} \A^-$.

The third and final possibility is that $su > u$ and $su \not\in W_c$,
meaning that $su$ is weakly complex.  Here, $\ell(su) = \ell(u) + 1
\leq \ell(w) < \ell(x)$, and by induction, $\te_{su} \in \latl{s}$.  
We therefore have $$
\te_{su} = \sum_{u' \in W_c} a'_{u'} \te_{u'}
,$$ where $su' < u'$ whenever $a'_{u'} \not\in v^{-1}\A^-$.  Since $sw > w$,
it follows that 
$\lrt{\te_w}{\te_{u'}} \in v^{-1} \A^-$.  By bilinearity, we have
$\lrt{\te_w}{\te_{su}} \in v^{-1} \A^-$.

We have now shown that $\te_x \in \L$.  Running through the argument again
with this in mind, we see that $\lrt{\te_x}{\te_u} \in v^{-1} \A^-$ unless
$su < u$, which by Lemma \secd.1 shows that $\te_x \in \latl{s}$.
\qed\enddemo

An interesting question is whether one can replace ``weakly complex'' in
Lemma \secd.9 by ``complex''; see \S\secconc\  below for more details.

\proclaim{Proposition \secd.10}
Suppose that $s, t \in S$ are noncommuting generators of
the Coxeter group $W$, and that $\te_x \in \latl{u}$ whenever $x$ is weakly
complex, $ux \in W_c$ and $u \in S$.  Let $w \in W_c$.   Then we have:
\item{\rm (i)}{$$
\te_s \te_w \in \cases
v \latl{s} & \text{ if } sw < w, \cr
\latl{s} & \text{ if } sw > w;
\endcases
$$} \item{\rm (ii)}
{$\te_s \L \cap \L \subseteq \latl{s}$;}
\item{\rm (iii)}
{$\te_s \latl{t} \subseteq \latl{st}$.}
\item{\rm (iv)}
{if $a \in S$ does not commute with $t$ and $a \ne s$, then 
$\te_a \latl{st} \subseteq \latl{a}$.}
\endproclaim

\demo{Proof}
If $sw > w$, then either $sw \in W_c$, in which case $\te_{sw} \in \latl{s}$
by definition, or $sw$ is weakly complex, in which case $\te_{sw} \in 
\latl{s}$ by hypothesis.  If, on the other hand, 
$sw < w$, we have $$
\te_s \te_w = \te_{sw} + (v - v^{-1}) \te_w
.$$  Part (i) follows because $sw, w \in W_c$.

For (ii), let $x \in \L$, and write $$
x = \sum_{u \in W_c} a_u \te_u
,$$ where $a_u \in \A^-$.  It follows from the proof of (i) that if $su < u$, 
we must have $a_u \in v^{-1} \A^-$: otherwise, the coefficient of $\te_u$
in $\te_s x$ would fail to lie in $\A^-$.  The claims of (ii) now follow
from the statement of (i).

Part (iii) follows from (i) and the fact that $tw < w$ and $sw < w$ are
mutually exclusive conditions for $w \in W_c$.  (This is because if
$tw < w$ and $sw < w$ then $w$ has a reduced expression beginning with
an alternating sequence of $m(s, t)$ occurrences of $s$ and $t$.)

For (iv), let $x' \in \latl{st}$, and write $$
x' = \sum_{u \in W_c} a'_u \te_u
,$$ where $a'_u \in \A^-$.  If $a'_u \not\in v^{-1} \A^-$, then $u$ has
a reduced expression beginning with $st$.  This means that $u$ cannot also
have a reduced expression beginning with $a$, because the reduced expressions
of $u$ are commutation equivalent and the leftmost $t$ will be to the left 
of any occurrence of $a$ in a reduced expression for $u$.  This means that
$au > u$, and thus $\te_a \te_u \in \latl{a}$ by part (i).  If, on the
other hand, $a'_u \in v^{-1} \A^-$, we have $\te_a \te_u \in v \latl{a}$.
The proof now follows.
\qed\enddemo

\proclaim{Lemma \secd.11}
Let $W$ be a Coxeter group and let $I = \{s, t\}$ be a pair of 
noncommuting generators in $S$.
Suppose that whenever $x$ is weakly complex, $ux \in W_c$ and $u \in S$, 
we have $\te_x \in \latl{u}$.  Let $w = w_I w^I$ be such that $w^I \in W_c$.  
\item{\rm (i)}{If $sw_I < w_I$, then $\te_w \in \latl{s}$.}
\item{\rm (ii)}{If $w_I = w_{st}$, the longest element in $W_I$, then 
$\te_{w_I} \te_{w^I} \in v^{-1} \latl{s}$, and $$
\te_{w_I} \te_{w^I} + v^{-1} \te_{sw_I} \te_{w^I} + v^{-1} \te_{t w_I} 
\te_{w^I} \in v^{-2} \L
.$$}\endproclaim

\demo{Proof}
We first prove (i), where the statement is trivial if $w_I = 1$.
Assume this is not the case.
The element $w_I$ has a reduced expression ending in $u \in S$, and 
Proposition \secd.10 (i) shows that $\te_u \te_{w^I} \in \latl{u}$.  We can
then repeatedly left multiply by other elements $\te_s$, appealing to 
Proposition \secd.10 (iii) to complete the proof.

Part (ii) follows by combining part (i) and equation (1) of Lemma \secd.7.
\qed\enddemo

\proclaim{Proposition \secd.12}
If the Coxeter group $W$ has Property B and Property F, then $W$ has
Property W.
\endproclaim

\demo{Proof}
Let $x$ be weakly complex, and write $x = sw$, where $w \in W_c$ and 
$s \in S$.  Let
$w = w_1 w_2 w_3$ be a reduced expression as in Lemma \secd.5 (i).

The proof is by induction on $\ell(w)$.  Since Property F holds, either
(i) $w$ is a product of commuting generators (which is incompatible with $x$
being weakly complex), or (ii) $w = abw'$ (where $a, b \in S$ are 
noncommuting generators) is left reducible to an element $y = bw'$, or 
(iii) $w = w'ba$ (with $a, b$ as before) is right reducible to an element 
$y = w'b$.

Suppose we are in case (ii) and $s$ fails to commute with $a$.  Since all
reduced expressions of $w$ are commutation equivalent, we must have $a = t$
and the element $w_1$ commutes with both $s$ and $t$.  This implies that
$(sw)_I = w_{st}$, where $I = \{s, t\}$.  By Lemma \secd.11 (ii), this shows 
that $\te_{sw} \in v^{-1} \L^-$, as required.

Suppose now that we are in case (ii) and $s$ commutes with $a$, but does not
commute with $b$.  This forces $b = t$ and $sw$ has a reduced expression of
the form $a w_{st} x'$, where $x' = (w_{st} x')^I$.  By Lemma \secd.11 (ii), 
we have $$\eqalign{
\te_{sw} &= \te_a \te_{w_{st}} \te_{x'} \cr
&= v^{-1} \te_a (- \te_{s w_{st}} \te_{x'} - \te_{t w_{st}} \te_{x'} + z),\cr
}$$ where $z \in v^{-1} \L$.  Proposition \secd.10 (i), which is
applicable by Lemma 4.9, shows that $$\te_a z \in \L.$$  Since $s w_{st}$ has
a reduced expression beginning in $t$, Lemma \secd.11 (i) shows that
$\te_{s w_{st}}\te_{x'} \in \latl{t}$.  Because $a$ does not commute 
with $t$, Proposition \secd.10 (iii) now shows that 
$$\te_a (\te_{s w_{st}} \te_{x'}) \in \L.$$  The element $t w_{st}$ has a
reduced expression starting with $st$.  Lemma \secd.11 (i) shows that
$\te_{st w_{st}}\te_{x'} \in \latl{t}$, and then Proposition 
\secd.10 (iii) shows that $\te_{t w_{st}}\te_{x'} \in \latl{st}.$  
By Proposition \secd.10 (iv)
and the fact that $a$ does not commute with $t$, we have 
$$\te_a (\te_{s w_{st}} \te_{x'}) \in \L.$$  Combining these observations
shows that $\te_{sw} \in v^{-1} \L$.

We are now either in the situation of case (ii) but where $s$ commutes with
$a$ and $b$, or in the situation of case (iii).  Both possibilities mean that
$sw$ has a reduced expression of the form $abx'$ or of the form $x'ba$, where
$a$ and $b$ are noncommuting generators.  

Suppose that $sw = abx'$.
Since $\te_{bx'} \in \latl{b}$ by Lemma \secd.9,
Proposition \secd.10 (iii) shows that $\te_{abx'} \in \latl{a}$.  It will 
therefore be enough to show that
if $z \in W_c$ with $z' = az < z$, then $\lrt{{\te_{abx}}}
{\te_z} \in v^{-1}\A^-$.  We apply Lemma \secc.7 (i) to show that $$
\lrt{\te_a \te_b \te_{x'}}{\te_a \te_{z'}} 
= 
\lrt{\te_b \te_{x'}}{\te_b \te_a \te_{z'}} 
+ \lrt{\te_{x'}}{\te_a \te_{z'}} 
- \lrt{\te_b \te_{x'}}{\te_{z'}} 
.$$  In the other case, where $sw = x'ba$, a similar argument using Lemma
\secc.7 (ii) shows that for $z' = za < z$ we have $$
\lrt{\te_{x'} \te_b \te_a }{\te_{z'} \te_a} 
= 
\lrt{\te_{x'} \te_b}{\te_{z'} \te_a \te_b} 
+ \lrt{\te_{x'}}{\te_{z'} \te_a} 
- \lrt{\te_{x'} \te_b}{\te_{z'}} 
.$$

There are several possibilities to consider.

The first case is that $x' \not\in W_c$.  
If $sw = abx'$, then by Lemma \secd.5 (iii), $x'$ must be weakly complex.  
This also implies that $bx'$ is weakly complex,
so $\te_b \te_{x'}$ and $\te_{x'}$ lie in $v^{-1} \L$ by induction.
By Proposition \secd.10 (i) and (iii), we see that $\te_b \te_a \te_{z'}$,
$\te_a \te_{z'}$ and $\te_{z'}$ all lie in $\L$.  This means that
$\lrt{\te_a \te_b \te_{x'}}{\te_a \te_{z'}}$ can be written as a sum of
three terms, each of which lies in $v^{-1} \A^-$, as required.
The alternative situation where $sw = x'ba$ and $x' \not\in W_c$ may be 
treated similarly.

If $sw = abx'$, it is not possible for $x' \in W_c$ and $bx' \not\in W_c$,
because the fact that $s$ commutes with $a$ and $b$ means that $a$ and $b$
correspond to generators in the factor $w_1$ of Lemma \secd.5 (i).  However,
if $sw = x'ba$, it is possible for $x' \in W_c$ and $x'b \not\in W_c$.
In this case, we may argue as before except as regards 
the term $\lrt{\te_{x'}}{\te_{z'} \te_a} = \lrt{\te_{x'}}{\te_z}$.  Since
$za < z$, this term will lie in $v^{-1} \A^-$ unless $x'a < x'$, in
other words, if $x' = x''a$ reduced.  Since
$x''a \in W_c$ and $x''ab \not\in W_c$, Lemma \secd.5 (i) shows that
$x''ab$ has a reduced expression of the form $x'''w_{ab}$.  This is a
contradiction, because it shows that $x'b$ has a reduced expression ending in
$a$, and yet $x'ba > x'b$.

If $sw = abx'$, the only other possibility is that $x', bx' \in W_c$ and 
$abx' \not\in W_c$.
Arguing as in the previous paragraph, $abx'$ has a reduced expression
beginning with $w_{ab}$.  The analysis of this case is now the same as
when $s$ fails to commute with $a$, which was considered above using Lemma
\secd.11 (ii).

The only remaining case is where $sw = x'ba$, $x', x'b \in W_c$ and 
$x'ba \not\in W_c$.  This may be treated analogously.
\qed\enddemo

\head \sece.  Inductive computation of the $\mu(x, w)$ \endhead

If $x, w \in W$, the Kazhdan--Lusztig polynomial $P_{x, w}$
is a polynomial in $q$, and if $x \ne w$, it has degree at most 
$(\ell(w) - \ell(x) - 1)/2$.  (If $x = w$, we have $P_{x, w} = 1$, and if
$x \not\leq w$ in the Bruhat order, we have $P_{x, w} = 0$.)  The cases
where the maximum degree bound is achieved are of 
particular importance.  (This can only happen when $\ell(w)$ and 
$\ell(x)$ are unequal modulo 2.)  If $x \ne w$, we denote the
coefficient of $q^{(\ell(w) - \ell(x) - 1)/2}$ in $P_{x,w}$ by $\mu(x, w)$.
Clearly, $\mu(x, w)$ will be zero unless $x < w$ and 
$\ell(w)$ and $\ell(x)$ are unequal modulo 2.

When $x, y \in W_c$, there are analogues $M(x, y)$ of the integers $\mu(x, y)$ 
associated to the basis $\{ c_w : w \in W_c \}$ of $TL(X)$.  These are
important for our purposes for two reasons: first, it often happens that
$M(x, y) = \mu(x, y)$, and secondly, the $M(x, y)$ are typically
much easier to compute than the $\mu(x, y)$ in general.  The goal of this
section is to relate the $M(x, y)$ to the structure constants of the 
basis $\{c_w : w \in W_c\}$ and to establish agreement, in certain cases, 
between the $M(x, y)$ and the $\mu(x, y)$.

As we shall see, one reason Property W is important is that 
it allows the inductive computation of the $c$-basis.

\definition{Definition \sece.1}
Let $W$ be any Coxeter group and let $y, w \in W_c$.  Let us write $$
\eqalignno{
c_w & = \sum_{y \in W_c} p^*(y, w) \te_y & (2)
}$$ and $$\eqalignno{
\te_w &= \sum_{y \in W_c} \e_y \e_w q^*(y, w) c_y
, & (3)}$$ where $\e_z$ means $(-1)^{\ell(z)}$.  
If $w \not\in W_c$ or $y \not\in W_c$, we make the convention that
$p^*(y, w) = 0$.  If $y \not\in W_c$, we define $q^*(y, w) = 0$; if $y \in W_c$
but $w \not\in W_c$, the formula (3) still makes sense, and we define
$q^*(y, w)$ as usual.  We also define 
$p(y, w) := v^{\ell(w) - \ell(y)} p^*(y, w)$ 
and
$q(y, w) := v^{\ell(w) - \ell(y)} q^*(y, w)$.
We define $M(y, w)$ to be the (integer) coefficient of $v^{-1}$ in $p^*(y, w)$,
and we write $y \prec w$ to mean that $M(y, w) \ne 0$.
\enddefinition

\proclaim{Lemma \sece.2}
Let $W$ be an arbitrary Coxeter group, and let $w, x \in W_c$.
\item{\rm (i)}{We have $p^*(w, w) = q^*(w, w) = 1$.}
\item{\rm (ii)}{If $x \not< w$, we have $p^*(x, w) = q^*(x, w) = 0$.}
\item{\rm (iii)}{If $x < w$, then $p^*(x, w)$ and $q^*(x, w)$ are 
elements of $v^{-1} \A^-$.}
\item{\rm (iv)}{The set $\{ v^{\ell(w)} c_w : w \in W_c \}$ is a 
$\zed[q, q^{-1}]$-basis for $TL_q(X)$.}
\item{\rm (v)}{The Laurent polynomials $p(x, w)$ and $q(x, w)$ lie in
$\zed[q, q^{-1}]$.}
\item{\rm (vi)}{If $x \prec w$ then $\e_x = -\e_w$.}
\item{\rm (vii)}{The coefficient of $v^{-1}$ in $q^*(x, w)$ is $M(x, w)$.}
\endproclaim

\demo{Proof}
Parts (i), (ii) and (iii) are immediate consequences of Theorem \secb.1 (ii).

We can uniquely write $c_w$ (or any 
element of $TL(X)$ as $c_w = x_1 + x_2$, where $x_1 \in TL_q(X)$
and $x_2 \in v TL_q(X)$.  By Theorem \secb.1 (ii), we have $x_1, x_2
\in \L$, and furthermore, we have $\pi(x_i) = \te_w$ and $\pi(x_j) = 0$
for $\{i, j\} = \{1, 2\}$.

The ring homomorphism $\bar{\ }$ fixes the $\zed[q, q^{-1}]$-algebras
$\H_q(X)$ and $TL_q(X)$, so the fact that $\overline{c_w} = c_w$ shows that
$\overline{x_i} = x_i$ for $i \in \{1, 2\}$.  The uniqueness properties of
$c_w$ now show that $x_i = c_w$ and $x_j = 0$.  Since $v^{\ell(w)} \te_w
\in TL_q(X)$, we now see that  $v^{\ell(w)} c_w \in TL_q(X)$.  Part (iv)
follows from these observations.

Since $v^{\ell(w)} c_w \in TL_q(X)$, it follows that 
$v^{\ell(w)} q^*(x, w) \te_x = q(x, w) t_x \in TL_q(X)$, from which
statement (v) for the $q(x, w)$ follows.  It follows easily from the 
definitions that $$\eqalignno{
\sum_{z \in W_c} p^*(x, z) (\e_z \e_w q^*(z, w)) &= \delta_{x, w}
, & (4) \cr}$$ and thus that $$\eqalignno{
\sum_{z \in W_c} v^{\ell(z) - \ell(x)} p^*(x, z) (\e_z \e_w 
v^{\ell(w) - \ell(z)} q^*(z, w)) &= \delta_{x, w}
, & (5)}$$  in other words, that the matrices $(p(x, w))$ and
$(\e_x \e_w q(x, w))$ are also mutually inverse.  Statement (v) for the
$p(x, w)$ follows from this.

It follows from (v) that if $\e_x = \e_w$ then $p^*(x, w)$ lies in $\zed[q]$,
and if $\e_x = -\e_w$ then $p^*(x, w)$ lies in $v \zed[q]$.  If $M(x, w) 
\ne 0$, this shows that $p^*(x, w)$ lies in $v \zed[q]$, and (vi) follows.

Define $M'(x, w)$ to be the coefficient of $v^{-1}$ in $q^*(x, w)$.
Equating coefficients of $v^{-1}$ on each side of (4) and applying (i), (ii)
and (iii), we find that $$
\e_w \e_w M(x, w) + \e_x \e_w M'(x, w) = 0
.$$  If $M(x, w) = 0$, then $M'(x, w) = 0$ as required.  If not, (vi) shows
that $\e_x = -\e_w$ and again $M(x, w) = M'(x, w)$, completing the proof.
\qed\enddemo

The following formulae are analogues of \cite{{\bf 18}, 1.0.a} and \cite{{\bf 20}, 4.3.1}.

\proclaim{Proposition \sece.3}
Suppose that the Coxeter group $W$ has Property W.  Let $w \in W_c$ and
$s \in S$.  Then we have $$
c_s c_w = \cases
(v + v^{-1}) c_w & \text{ if } \ell(sw) < \ell(w),\cr
c_{sw} + \sum_{sy < y} M(y, w) c_y
& \text{ if } \ell(sw) > \ell(w),\cr
\endcases$$ where $c_z$ is defined to be zero whenever $z \not\in W_c$.
\endproclaim

\demo{Proof}
Let us observe that the basis element $c_1$ is the identity element of $TL(X)$,
and that if $s \in S$, we have $c_s = v^{-1} \te_1 + \te_s$.  These claims
can be proved by checking the uniqueness criteria of Theorem \secb.1 (ii).

We first deal with the case where $sw > w$.  From Theorem \secb.1 and
Definition \sece.1, we know that $$
c_w = \te_w + \sum_{{y < w} \atop {y \in W_c}} p^*(y, w) \te_y
,$$ where the coefficient of $v^{-1}$ in $p^*(y, w)$ is $M(y, w)$.  It follows
that $$
\te_s c_w = (\te_s \te_w) + \sum_{{y < w} \atop {y \in W_c}} p^*(y, w) 
(\te_s \te_y)
.$$  Proposition \secd.10 (i) and the fact that the $p^*(y, w)$ lie in 
$v^{-1} \A^-$ show that $\te_s c_w \in \L$.  

Since $c_s = v^{-1} \te_1 + \te_s$, we have
$c_s c_w \in \L$.  Since $\bar{\ }$ is a ring homomorphism, Theorem \secb.1 (ii)
shows that $\overline{c_s c_w} = c_s c_w$, and Theorem \secb.1 (iv) shows
that it is enough to prove that $$
\pi(c_s c_w) = \pi\left( c_{sw} + \sum_{sy < y} M(y, w) c_y\right)
.$$  Using the above formula for $c_s$, this is equivalent to $$
\pi(\te_s c_w) = \pi\left( c_{sw} + \sum_{sy < y} M(y, w) \te_y\right)
.$$  If $sw \not\in W_c$, then $c_{sw}$ is defined to be zero, and 
$\pi(\te_s \te_w) = 0$ by Property W.  If, on the other hand, 
$sw \in W_c$, we have $\pi(c_{sw}) = \pi(\te_s \te_w)$ by Theorem \secb.1 (ii).
Suppose that $y < w$.
If $sy > y$, we have $\te_s \te_y \in \L$ by Proposition \secd.10 (i), and thus
$\pi(p^*(y, w) \te_s \te_y) = 0$.  If, on the other hand, $sy < y$, we have
$\te_s \te_y = (v - v^{-1}) \te_y + \te_{sy}$, which implies that $$
\pi(\te_s \te_y) = \pi((v - v^{-1})p^*(y, w) \te_y) = \pi(M(y, w) \te_y)
.$$  The result now follows from the formula for $\te_s c_w$.

It remains to show that $c_s c_w = (v + v^{-1}) c_w$ if $sw < w$, which we will
prove by induction on $\ell(w)$.  The case $\ell(w) = 0$ cannot occur,
and the case $\ell(w) = 1$ follows from the Hecke algebra identity $$
C'_s C'_s = (v + v^{-1}) C'_s
.$$  Suppose now that $\ell(w) > 1$, and write $w = sx$.  We now know that $$
c_{sx} = c_s c_x - \sum_{sy < y} M(y, x) c_y
.$$  Since $y < x$ for each $y$ appearing in the sum with nonzero coefficient,
we have $c_s c_y = (v + v^{-1}) c_y$ by induction.  We also have $c_s c_s c_x
= (v + v^{-1}) c_s c_x$ by induction, from which the claim follows.
\qed\enddemo

\proclaim{Corollary \sece.4}
Suppose that $W$ has Property W.  Then the set $$
\{ x \in TL(X) : c_s x = (v + v^{-1}) x \}
$$ is the free $\A$-submodule of $TL(X)$ with basis $\{c_y : sy < y\}$.
\endproclaim

\demo{Proof}
This is immediate from Proposition \sece.3 and the observation that 
all the basis elements $c_y$ appearing in the expression for $c_s c_w$
in that result satisfy $sy < y$.
\qed\enddemo

\proclaim{Lemma \sece.5}
Suppose that $W$ has Property W.  Let $x, w \in W_c$ and $s
\in S$ be such that $sw > w$ (although we do not assume that
$sw \in W_c$).
\item{\rm (i)}
{If $sx > x$ then we have $q(x, sw) = q(x, w)$.}
\item{\rm (ii)}
{If $sx < x$ then we have $$\eqalignno{
q(x, sw) &= -v^2 q(x, w) + q(sx, w) +
 \sum_{{x \prec y \leq w} \atop {sy > y}} v^{\ell(y) + 1 - \ell(x)}
M(x, y) q(y, w)
. & (6)}$$}
\endproclaim

\demo{Proof}
Using (3) and Proposition \sece.3 we find that $$\eqalign{
v^{-1} \te_w + \te_{sw} =& c_s \te_w \cr
=& \sum_{x \leq w} \e_x \e_w q^*(x, w) c_s c_x\cr
=& \left( \sum_{{x \leq w} \atop {sx < x}} 
\e_x \e_w (v + v^{-1}) q^*(x, w) c_x \right) \cr
&+ \sum_{{x \leq w} \atop {sx > x}} \e_x \e_w q^*(x, w) 
\left( 
c_{sx} + \sum_{{z \prec x} \atop {sz < z}} M(z, w) c_z
\right).\cr
}$$  Using (3) again to equate the coefficients of $c_x$ on each side of the
equation, routine calculations yield the stated identities.
\qed\enddemo

\proclaim{Lemma \sece.6}
Suppose that $W$ has Property W, and let $x, w \in W_c$.
\item{\rm (i)}
{The $q(x, w)$ and $p(x, w)$ are polynomials in $q$, and $q(x, w)$ has 
constant term $1$.}
\item{\rm (ii)}
{If $x < w$, the $q(x, w)$ and $p(x, w)$ have degree at most 
$(\ell(w) - \ell(x) - 1)/2$ as polynomials in $q$, with the degree bound
being attained if and only if $M(x, w) \ne 0$.}
\item{\rm (iii)}
{Let $x, w \in W_c$ and $s \in S$ be such that $sw < w$ and $sx > x$.
If $M(x, w) \ne 0$ then we must have $x = sw$ and $M(x, w) = 1$.}
\endproclaim

\demo{Proof}
We prove (i) by induction on $\ell(w)$.  The case $\ell(w) = 0$ follows
from Lemma \sece.2 (i).  For the inductive step, we write $w = sw'$ for
some $s \in S$ with $w' < w$.  The assertions of (i) for the $q(x, w)$
follow quickly from the observation that the quantity $\ell(y) + 1 - \ell(x)$ 
appearing in the sum of Lemma \sece.5 (ii) is a strictly positive even integer.
The assertions about the $p(x, w)$ then follow from equation (5), Lemma
\sece.2 (i) and linear algebra.

Part (ii) follows from (i) and Lemma \sece.2 (iii).

For (iii), Lemma \sece.5 (i) shows that $q(x, w') = q(x, w)$.  If $x \prec w$,
so that $x < w$, then the degree of $q(x, w)$ must be 
$(\ell(w) - \ell(x) - 1)/2$ by (ii).  This exceeds the degree bound of
$(\ell(w') - \ell(x) - 1)/2$ which would apply to $q(x, w')$ 
unless $x = w'$, as required.
\qed\enddemo

\remark{Remark \sece.7}
Unlike the case of the polynomials $P_{x, w}$, it is not true that
$p(x, w)$ has constant term $1$.  If this were the case, equation (5) and 
the argument of
\cite{{\bf 16}, Corollary 7.13} would show that for $x, w \in W_c$, each interval $$
\{ y \in W_c: x \leq y \leq w \}
$$ would contain equal numbers of elements of odd and even lengths.  However,
this is not true in type $A_3$: take $x = s_2$ and $w = s_2 s_1 s_3 s_2$.
\endremark

\definition{Definition \sece.8}
As in \cite{{\bf 20}}, we define $$
\tmu(x, y) = \cases
\mu(x, y) & \text{ if } x \leq y;\cr
\mu(y, x) & \text{ if } x > y.\cr
\endcases$$  Analogously, we define $$
\tem(x, y) = \cases
M(x, y) & \text{ if } x \leq y;\cr
M(y, x) & \text{ if } x > y.\cr
\endcases$$
\enddefinition

In order to show that the coefficients $M(x, y)$ appearing in Lemma \sece.6 
are equal to the coefficients $\mu(x, y)$ of \cite{{\bf 18}}, 
we show that each set of
coefficients satisfies a common recurrence relation.  This recurrence relation
is easy to explain in terms of star operations.

\proclaim{Proposition \sece.9 (Lusztig)}
Let $W$ be an arbitrary Coxeter group, and let $x$ and $w$ be 
elements of $\{s, t\}$-strings (for the same $s$ and $t$, but 
possibly different strings).  Suppose that 
$\ldescent{x} \cap \{s, t\} \ne \ldescent{w} \cap \{s, t\}$.
Then $$
\tmu(_*x, w) + \tmu(^*x, w) = \tmu(x, {_*w}) + \tmu(x, {^*w})
,$$ where we define $\tmu(a, b) = 0$ if either $a$ or $b$ is an undefined
symbol.
\endproclaim

\demo{Proof}
This result is implicit in \cite{{\bf 20}, \S10.4}, and is what Lusztig is referring
to by ``an analogous result holds for arbitrary $m$''.  (A proof may also be
obtained by modifying the argument below (Proposition \sece.12)
for the symbols $\tem(x, y)$.)
\qed\enddemo

The following is a routine exercise using the subexpression characterization
of the Bruhat order of a Coxeter group (see also \cite{{\bf 1}, Proposition 2.5.1}).

\proclaim{Lemma \sece.10}
Let $W$ be an arbitrary Coxeter group, let $I$ be as in Definition \secc.1
and let $x = x_I x^I$, $y = y_I y^I$, $w = w_I w^I$ be three elements of
$W$.  If $x \leq w$ then we must have $x^I \leq w^I$.
Furthermore, if $x^I = w^I$ and $x \leq y \leq w$, we must have $x^I = y^I
= w^I$ and $x_I \leq y_I \leq w_I$.
\qed\endproclaim

\proclaim{Lemma \sece.11}
Suppose that $W$ satisfies Property W.
Let $x = x_I w^I$ and $w = w_I w^I$ be two elements of $W_c$ in the same
coset of $W_I$, where $I$ is as in Definition \secc.1.  Then we have 
$q(x, w) = q(x_I, w_I).$
\endproclaim

\demo{Proof}
By Lemma \sece.2 (ii), we may assume $x \leq w$, which implies $x_I \leq w_I$
by Lemma \sece.10.  We will
proceed by induction on $\ell(w_I)$.  If $\ell(w_I) = 0$ then necessarily
$x = w$ and $x_I = w_I$, and the statement follows from Lemma \sece.2 (i).
If $\ell(w_I) > 0$, write $w_I = sw_I' > w_I' \in W_c$, where $s \in I$.  This 
implies that $w' = sw < w$.  

Suppose that $sx_I > x_I$; this implies that
$sx > x$.  Lemma \sece.5 now shows that $$
q(x, w) = q(x_I, w'_I) = q(x_I, sw'_I) = q(x_I, w_I)
,$$ by induction.

Now suppose that $sx_I < x_I$, which means that $sx_I$ and $sx$ lie in $W_c$.
By equation (6) and Lemma \sece.10, we have $$\eqalign{
q(x, sw') 
&= -v^2 q(x, w') + q(sx, w') +
 \sum_{{x \prec y \leq w'} \atop {sy > y}} v^{\ell(y) + 1 - \ell(x)}
M(x, y) q(y, w') \cr
&= -v^2 q(x_I, w'_I) + q(sx_I, w'_I) +
 \sum_{{x_I \prec y_I \leq w'_I} \atop {sy_I > y_I}} v^{\ell(y) + 1 - \ell(x)}
M(x, y) q(y_I, w'_I) \cr
&= q(x_I, sw'_I),
}$$ as required.
\qed\enddemo

\proclaim{Proposition \sece.12}
Suppose that the Coxeter group $W$ satisfies Property F and Property W.
Let $x, w \in W_c$ be
elements of $\{s, t\}$-strings (for the same $s$ and $t$, but possibly 
different strings) and let $I = \{s, t\}$.  
Suppose that 
$\ldescent{x} \cap \{s, t\} \ne \ldescent{w} \cap \{s, t\}$.  Then $$
\tem({_*x}, w) + \tem(^*x, w) = \tem(x, {_*w}) + \tem(x, {^*w})
,$$ where we define $\tem(a, b) = 0$ if either $a$ or $b$ is an undefined
symbol.  Furthermore, if $x^I \ne w^I$ and $\ell(x) \leq \ell(w)$, 
we can replace $\tem(a, b)$ by $M(a, b)$ throughout.
\endproclaim

\demo{Proof}
Note that the elements ${_*x}$, $x$ and ${^*x}$ have the same coset 
representative, $x^I$, and that the elements ${_*w}$, $w$ and ${^*w}$ have
the same coset representative, $w^I$.

We may assume that $\e_x = \e_w$ throughout, otherwise all terms are zero by 
Lemma \sece.2 (vi).

Suppose first that $x^I$ = $w^I$.  By Lemma \sece.11, it is enough to verify
the statement when $x$ and $y$ are replaced by $x_I$ and $w_I$, respectively;
in other words, $W$ may be assumed to be a dihedral group.  In this case
it is easily checked that the unique solution to the identities in Lemma
\sece.5 is $$
q(x, w) = \cases
1 & \text{\ if\ } x \leq w;\cr
0 & \text{\ otherwise.}\cr
\endcases$$  We therefore have, for $a, b \in W_I \cap W_c$, 
$\tem(a, b) = 1$ if and only if $\ell(b) = \ell(a) \pm 1$.
Verification of the claim is now an easy case by case check
according to the value of $\ell(x) - \ell(w)$.

Now suppose that $x^I \ne w^I$.  To fix notation, let us suppose that
$sw < w$, and thus $tx < x$.  By Lemma \sece.5 (ii), we have $$\eqalignno{
q(x, tw) &= -v^2 q(x, w) + q(tx, w) +
 \sum_{{x \prec y \leq w} \atop {ty > y}} v^{\ell(y) + 1 - \ell(x)}
M(x, y) q(y, w)
. & (7)}$$  By Lemma \sece.5, we may replace $q(x, w)$ in equation (7) with
$q(x, sw)$, which expresses (7) as a sum of terms each of which is a 
polynomial in $q$ of degree at most $(\ell(w) - \ell(x))/2$.  
Suppose first that $tw \not\in W_c$, in other words, that $tw$ has a
reduced expression beginning with $w_{st}$ and that ${^*w}$ is not defined.  
After the substitution just described, (7) shows 
that $q(x, tw)$ has degree at most $(\ell(tw) - \ell(x) - 1)/2$ as
a polynomial in $q$.  If this degree bound is attained, we find that
$v^{-1}$ appears with nonzero coefficient in $q^*(x, tw)$.  Lemma \secd.11 (ii)
shows that this can only happen if either $x = w$, or if $x = stw$ and
$stw \in W_c$.  However, both these possibilities imply that $x^I = w^I$,
and this case has already been eliminated.

We may now assume that $tw = {^*w}$, and hence that $tw \in W_c$.  Considering
the coefficients of $q^{(\ell(w) - \ell(x))/2}$ in (7), we find that 
$$\eqalignno{
M(x, {^*w}) &= - M(x, sw) + M(tx, w) + \sum_{{x \prec y \leq w} \atop {ty > y}}
M(x, y) M(y, w)
. & (8)}$$  Suppose that $M(x, y) M(y, w)$ is a nonzero term in the sum of
equation (8).  We know that $s \in \ldescent{w}$.  By Lemma \sece.6 (iii),
this means that either $y = sw$, or that $s \in \ldescent{y}$.  In the
latter case, we can apply Lemma \sece.6 (iii) again to see that either
$s \in \ldescent{x}$ or $x = sy$.  However, we have seen that $tx < x$,
and since $x$ lies in an $\{s, t\}$-string, this forces $sx > x$.  There 
are thus only two possibilies for values of $y$ giving nonzero terms in the 
sum, namely $y = sw$ or $y = sx$.

Consider first the case where $y = sx$.  Since $ty > y$ for all $y$ in
the sum, we have $tsx > x$.  Since $x < sx < tsx$, this means that
$sx = {^*x}$.  In any case, we have a contribution of $M({^*x}, w)$ to
the sum in (8).

Now consider the case where $y = sw$.  As above, we have $ty > y$ and
thus $tsw > sw$.  We have observed that $sw < w$, and this means that
$sw$ is not an element of the $\{s, t\}$-string containing $w$, or equivalently
that ${_*w}$ is not defined.  The term
$y = sw$ contributes a term $M(x, sw)$ to the sum, and this cancels the
term $-M(x, sw)$ already appearing.  This produces a total of
$-M(x, {_*w})$, \idest zero.  

On the other hand, if ${_*w}$ is defined, we must have
${_*w} = sw$ and $tsw < sw$.  This means that the case $y = sw$ cannot occur, 
and the term $-M(x, sw) = -M(x, {_*w})$ already appearing in (8) is not
cancelled by a term in the sum, again leaving a total contribution of 
$-M(x, {_*w})$.

It remains to consider the term $M(tx, w)$ appearing in (8).  We know that
$sw < w$, so for $M(tx, w) \ne 0$, we require either $tx = sw$, or $stx < tx$.
If $tx = sw$ then $x^I = w^I$, and we have already eliminated this case.
If, on the other hand, $stx < tx$, then we have $tx = {_*x}$.  In any case,
we find that $M(tx, w) = M({_*x}, w)$.

In summary, we have transformed (8) into the equation $$
M(x, {^*w}) = -M(x, {_*w}) + M({_*x}, w) + M({^*x}, w)
,$$ from which the claims follow.
\qed\enddemo

\proclaim{Theorem \sece.13}
Suppose that the Coxeter group $W$ satisfies Property F and Property W,
and let $x, w \in W_c$.  Then $M(x, w) = \mu(x, w)$, and in particular, 
we have $$
c_s c_w = \cases
(v + v^{-1}) c_w & \text{ if } \ell(sw) < \ell(w),\cr
c_{sw} + \sum_{sy < y} \mu(y, w) c_y
& \text{ if } \ell(sw) > \ell(w),\cr
\endcases$$ where $c_z$ is defined to be zero whenever $z \not\in W_c$.
\endproclaim

\demo{Proof}
The second claim is immediate from the first and Proposition \sece.3.

Let us first consider the case where $w = s_1 s_2 \cdots s_r$ is a product of 
distinct commuting generators.  In this case, direct computation shows that $$
C'_w = C'_{s_1} C'_{s_2} \cdots C'_{s_r}
$$ and $$
c_w = c_{s_1} c_{s_2} \cdots c_{s_r}
,$$ from which it follows (by considering the coefficient of $\T_w$ or $\te_w$
on the right hand sides of the equations) that $$
M(x, w) = \mu(x, w) = \cases
1 & \text{ if } x < w \text{ and } \ell(x) = \ell(w) - 1;\cr
0 & \text{ otherwise.}\cr
\endcases$$

We complete the proof of the first claim for $\mu(x, w)$ by induction on 
$\ell(w) - \ell(x)$.
The claim is trivial unless $\ell(w) - \ell(x)$ is an odd positive integer,
by Lemma \sece.2 (ii), (vi) and \cite{{\bf 18}, Definition 1.2}.  If $\ell(w)
= \ell(x) + 1$, Lemma \sece.2 (ii) shows that $M(x, w) = 0$ if $x \not< w$,
and Lemma \sece.6 (i) shows that $M(x, w) = 1$ if $x < w$.  The same is
true of the $\mu(x, w)$ by \cite{{\bf 18}, Definition 1.2, Lemma 2.6 (i)}.

For the inductive step, we may assume that $\ell(w) - \ell(x) > 3$.  Since
Property F holds and we have dealt with the case where $w$ is a product of
commuting generators, we may write $w = stw'$ or $w = w'ts$ reduced, where $s$
and $t$ are noncommuting generators.  We treat the former case; the latter
is dealt with by a symmetrical argument.  Since $w \in W_c$, we have
$w = {^*y}$, where $y = tw'$.  It suffices to compute $\tem(x, {^*y})$.
If $\ldescent{{^*y}} \not\subseteq \ldescent{x}$, Lemma \sece.5 (i)
shows that either $\tem(x, {^*y}) = 0$ or $\ell(x) = \ell({^*y}) - 1$, and
the latter case has already been dealt with.  Since $sw < w$, we may now
assume that $sx < x$, and since $x \in W_c$, we must have $tx > x$.  The
hypotheses of Proposition \sece.12 are now satisfied, and we use the
relation there to compute $\tem(x, {^*y})$ by induction.
The $\mu(x, w)$ satisfy the same recurrence, except that one uses 
\cite{{\bf 18}, (2.3e)} in place of Lemma \sece.5 (i), and Proposition
\sece.9 in place of Proposition \sece.12.
\qed\enddemo

\remark{Remark \sece.14}
Theorem \sece.13 was first observed in the $ADE$ case by Graham 
\cite{{\bf 7}, Theorem 9.9}, prior to the definition of the $c_w$-basis \cite{{\bf 13}}.
\endremark

\head \secf.  Positivity properties for the $c$-basis \endhead

In this section, we show how Property F and Property W may be used 
prove the positivity of structure constants for the $c$-basis, a property known to hold in all cases where the $c$-basis has been explicitly constructed.  
If Property S also holds, this gives an elementary proof that certain of 
the structure constants for the Kazhdan--Lusztig basis are positive.

The following well-known consequence of \cite{{\bf 18}, Theorem 1.3} 
is the model for Theorem \sece.13.

\proclaim{Lemma \secf.1 (Kazhdan--Lusztig)}
If $W$ is an arbitrary Coxeter group, then we have $$
C'_s C'_w = \cases
(v + v^{-1}) C'_w & \text{ if } sw < w;\cr
C'_{sw} + \sum_{{z \prec w} \atop {sz < z}} \mu(z, w) C'_z & \text { if }
sw > w.\cr
\endcases
$$  
\qed\endproclaim

\proclaim{Lemma \secf.2}
Let $W$ be an arbitrary Coxeter group, and let $I = \{s, t\} \in S$ be 
noncommuting generators and $w \in W_c$ be such that $tw < w$ and $sw > w$.  
Then we have $$
C'_s C'_w = C'_{sw} + C'_{_*w} + \sum_{I \subseteq \ldescent{z}} 
\mu(z, w) C'_z
,$$ where we interpret $C'_z$ to mean zero if $z$ is an undefined symbol.
In particular, we have $$\eqalignno{
C'_s C'_w &= C'_{^*w} + C'_{_*w} \mod J(X). & (9)\cr
}$$  
\endproclaim

\demo{Proof}
We use the formula of Lemma \secf.1 in the case where $sw > w$.
Now $tw < w$, so in order to have $z \prec w$, \cite{{\bf 18}, (2.3e)} shows that 
we need either $tz < z$ or 
$z = tw$.  If $tz < z$ then $z$ satisfies the conditions of the sum in 
the statement.  If $z = tw < w$ then $tz > z$ and $sz < z$, so
$z = {_*w}$, and $\mu(z, w) = 1$ by \cite{{\bf 18}, (2.3e)}.  The first
assertion now follows.

Suppose that $x \in W$ is such that $sx < x$ and $tx < x$.  Since
$\T_u C'_x = v C'_x$ for $u \in I$, an inductive argument using the formula
for $C'_{w_{st}}$ in terms of the $\T$-basis shows that $$
C'_{w_{st}} C'_x = (v + v^{-1}) (v^{m-1} + v^{m-3} + \cdots
+ v^{-(m-1)}) C'_x
,$$ where $m$ is the order of $st$.  (Note that if $m$ is infinite,
the hypotheses $sx < x$ and $tx < x$ are incompatible.)  Since $TL(X)$ is a
free $\A$-module, this shows
that $C'_x \in J(X)$.  Similarly, if ${^*w}$ is not defined, $C'_{^*w}
\in J(X)$.  The second assertion now follows.
\qed\enddemo

\proclaim{Proposition \secf.3}
Suppose that the Coxeter group $W$ satisfies Property F and Property W.
\item{\rm (i)}{The map 
$$\th : \H(X) \ra TL(X)$$
satisfies $\th(C'_w) = c_w$ whenever $w \in W_c$.}
\item{\rm (ii)}{If $I = \{s, t\}$ is a pair of noncommuting generators,
and we have $w \in W_c$ with $tw < w$, then we have $$
c_s c_w = c_{^*w} + c_{_*w}
.$$}
\endproclaim

\demo{Proof}
The proof of (i) is by induction on the length of $w$, the base case being 
where $w$ is a product of commuting generators.  If this is the case, and
$w = s_1 s_2 \cdots s_r$, it may be checked directly that $$
C'_w = \sum_{z < w} v^{\ell(z) - \ell(w)} \T_z
,$$ and because all the $z < w$ in the sum satisfy $z \in W_c$, it follows
that $$
c_w = \sum_{z < w} v^{\ell(z) - \ell(w)} \te_z
,$$ \idest $\th(C'_w) = c_w$.

Suppose that $w$ is not a product of commuting generators.  By Property F,
$w$ is either left star reducible or right star reducible.  We treat only
the case of left star reducibility, as the other is similar.

In this case, we can write $w = sx$ reduced, where $x \in W_c$ and $tx < x$ 
for some noncommuting generators $s$ and $t$.  Lemma \secf.2 shows that $$
C'_s C'_x = C'_{^*x} + C'_{_*x} \mod J(X)
.$$  Applying Theorem \sece.13, we find that $$
c_s c_x = c_{^*x} + c_{_*x}:$$
the reason for this is that the conditions $\mu(y, x) \ne 0$, $y \in W_c$, 
$tx < x$ and $sy < y$ force $ty > y$, $y = tx$ and $\mu(y, w) = 1$ by
Lemma \sece.6 (iii).  This completes the induction and the proof of (i).

Part (ii) follows from (i) and Lemma \secf.2.
\qed\enddemo


In order to prove positivity of structure constants, it is necessary to
have a good understanding of what happens in the much simpler case of
dihedral groups.  Let $I = \{s, t\}$ and let $W$ be the group of type
$I_2(m)$ generated by $I$.  
We define the Chebyshev polynomials of the second kind
to be the elements of $\zed[x]$ given by the conditions $P_0(x) = 1$,
$P_1(x) = x$ and $$\eqalignno{
P_n(x) &= xP_{n-1}(x) - P_{n-2}(x) & (10)\cr
}$$ for $n \geq 2$.  If $f(x) \in \zed[x]$, we define $f^{s, t}(x)$ to be 
the element of $\H$ given by the linear extension of the map sending
$x^n$ to the product $$
\underbrace{C'_s C'_t C'_s \ldots}_{n \text{ factors}}
$$ of alternating factors starting with $C'_s$.

\proclaim{Lemma \secf.4}
Let $W$ be a Coxeter group of type $I_2(m)$, and maintain the above
notation.  Then the $C'$-basis of $\H$ is given by the set $$\eqalign{
\{1\} 
& \cup \{(xP_i)^{s, t}(x) : i = 0, 1, \ldots, m-2\}\cr
& \cup \{(xP_i)^{t, s}(x) : i = 0, 1, \ldots, m-2\}\cr
& \cup \{(xP_{m-1})^{s, t}(x) = (xP_{m-1})^{t, s}(x)\}.\cr
}$$
\endproclaim

\demo{Proof}
This follows by a routine induction on $\ell(w)$ using Lemma \secf.1, equation
(10), and the fact that in type $I_2(m)$, we have 
$\mu(y, w) = 1$ if $\ell(y) = \ell(w) - 1$ and
$\mu(y, w) = 0$ otherwise. 
\qed\enddemo

\proclaim{Corollary \secf.5}
If $W$ is a Coxeter group of type $I_2(m)$, the $c$-basis of $TL(X)$
is given by the images under $\th$ of $$
\{1\} 
\cup \{(xP_i)^{s, t}(x) : i = 0, 1, \ldots, m-2\}
\cup \{(xP_i)^{t, s}(x) : i = 0, 1, \ldots, m-2\}
.$$
\endproclaim

\demo{Proof}
In this case, the ideal $J(X)$ is spanned by $$
C'_{w_0} = 
(xP_{m-1})^{s, t}(x) = (xP_{m-1})^{t, s}(x)
,$$ and the result now follows.
\qed\enddemo

The following result, which establishes positivity of structure constants
in the easy case of $TL(I_2(m))$, is our basic tool for proving positivity
in general.  Since the Laurent polynomial $v + v^{-1}$ appears frequently,
we will denote it by $\d$ from now on.

\proclaim{Proposition \secf.6}
Let $W$ be a Coxeter group of type $I_2(m)$, let $a, b \in W_c$ and write $$
c_a c_b = \sum_{w \in W_c} \l_w c_w
.$$  
\item{\rm (i)}{We have $\l_w \in \zed^{\geq 0}$ if $\rdescent{a} \cap 
\ldescent{b}
= \emptyset$, and $\l_w \in \d\zed^{\geq 0}$ otherwise.}
\item{\rm (ii)}{If $a \ne 1$, $b \ne 1$ and $\l_w \ne 0$, we have 
$\ldescent{w} = \ldescent{a}$ and $\rdescent{w} = \rdescent{b}$.}
\endproclaim

\demo{Proof}
If $a = 1$ or $b = 1$, the claims are clear, so suppose that this is not the
case.

Let $0 \leq i, j < m-1$, and let $K$ be the ideal $\lan P_{m-1}(x) \ran$ of 
$\zed[x]$.  If we write $$
P_i(x) P_j(x) = \sum_{0 \leq k < m-1} f_{i, j}^k P_k(x) \mod K
,$$ then it is well known (see, for example, \cite{{\bf 10}, Proposition 1.2.3}) 
that the $f_{i, j}^k$ lie in $\zed^{\geq 0}$, and furthermore, that 
$f_{i, j}^k \ne 0$ implies that $k \equiv i + j \mod 2$.

Because $x = P_1(x)$, we also see that $P_i(x) x P_j(x)$ can be written
as a positive combination of elements $P_k(x) \mod K$, and thus that
$x P_i(x) x P_j(x)$ can be written as a linear combination of 
$x P_k(x) \mod K$.  
The case in (i) where $\rdescent{a} \cap \ldescent{b} = \emptyset$ 
follows from this, and the ideal $K$ corresponds to the ideal $J(X)$.

It also follows that the product $(x P_i(x))P_j(x)$ can be written
as a positive combination of elements $x P_k(x) \mod K$.  
The other case of (i) follows from this observation.

The claims of (ii) follow by applying the fact that $k \equiv i + j \mod 2$
from above to the $c$-basis.
\qed\enddemo

The following result provides a convenient recursive method for computing
the $c$-basis.  

\proclaim{Lemma \secf.7}
Suppose $W$ satisfies Property F and Property W, and let $I = \{s, t\}$ be
a pair of noncommuting generators.  Let $w \in W_c$, let $w_I w^I$ be the
coset decomposition of $w$, and let $u \in I$ be the unique element of
$\rdescent{w_I}$.  Then we have $$
c_{w_I} c_{u w^I} = \d c_w
.$$
\endproclaim

\demo{Proof}
By Corollary \secf.5, we have an explicit expression for $c_{w_I}$, and
by Theorem \sece.13, we know that $c_u c_{u w^I} = \d c_{u w_I}$.
The proof follows by induction on $\ell(w_I)$, by applying Proposition \secf.3
(i) to equation (9), and comparing with equation (10).
\qed\enddemo


\proclaim{Lemma \secf.8}
Suppose $W$ satisfies Property F and Property W, and let $I = \{s, t\}$ be
a pair of noncommuting generators.  Let $1 \ne x \in W_c \cap W_I$ and 
$y \in W_c$ be such that $\rdescent{x} \subseteq \ldescent{y}$.  Writing $$
c_x c_y = \sum_{y \in W_c} f(x, y, w) c_w
,$$ we have $f(x, y, w) \in \d \zed^{\geq 0}$ for all $w$.
\endproclaim

\demo{Proof}
Let us write $y = y_I y^I$ and $u \in \rdescent{y_I}$, as in Lemma \secf.7.
Applying Lemma \secf.7, we see that $$
c_y = \d^{-1} c_{y_I} c_{uy^I} 
,$$ and thus $$
c_x c_y = (\d^{-1} c_x c_{y_I}) c_{uy^I}
.$$

The hypotheses of the statement require that $\rdescent{x} \cap \ldescent{y_I}
\ne \emptyset$, so by Proposition \secf.6 we have $$
\d^{-1} c_x c_{y_I} = \sum_{z \in W_c \cap W_I} \l_z c_z
,$$ where $\l_z \in \zed^{\geq 0}$ and $\l_z \ne 0$ implies that 
$\rdescent{c_z} = \{u\}$.  We can now apply Lemma \secf.7 to each term $c_z$
where $\l_z \ne 0$ to obtain $$
c_z c_{uy^I} = \d c_{zy^I}
.$$  Lemma \secd.5 (i) together with the fact that $zy^I$ is reduced shows that
$zy^I \in W_c$.  Putting all this together, we find that $$
c_x c_y = \sum_{z \in W_c \cap W_I} \d \l_z c_{zy^I}
,$$ which proves the statement.
\qed\enddemo


The next step is to show that the integers $\mu(y, w)$ appearing in the 
statement of Theorem \sece.13 are positive.  This is not obvious from the 
recurrence relations of propositions \sece.9 and \sece.12, except in easy cases
such as when the Coxeter graph is simply laced.  Note also that the $\mu(y, w)$
we are considering are not arbitrary: the set $\ldescent{y}$ properly
contains the set $\ldescent{w}$.

\proclaim{Lemma \secf.9}
Suppose that $W$ has Property F, and let $w \in W_c$ and 
$x = sw > w$.
Then one of the following situations must occur:
\item{\rm (i)}{$x$ is a product of commuting generators;}
\item{\rm (ii)}{$x \in W_c$ and there exists $I = \{s, t\} \subseteq S$ 
with $st \ne ts$ such that when $x = x_I x^I$, we have $\ell(x_I) > 1$;}
\item{\rm (iii)}{$x$ is weakly complex and has a reduced expression begining
with $w_{st}$ for some $t \in S$ with $st \ne ts$;}
\item{\rm (iv)}{there exists $I = \{u, u'\} \subset S$ with $s \not\in I$,
$uu' \ne u'u$, $su = us$ and $su' = u's$ such that when we write 
$w = w_I w^I$, we have $\ell(w_I) > 1$;}
\item{\rm (v)}{there exists $I = \{u, u'\} \subset S$ with 
$uu' \ne u'u$ such that when we write 
$w = (^Iw)(_Iw)$, we have $\ell(_Iw) > 1$;}
\item{\rm (vi)}{$x$ is weakly complex and there exist $t, u \in S$ with
$st \ne ts$, $ut \ne tu$ and $su = us$ such that $w$ has a reduced
expression of the form $$
u (tsts\cdots) x'
,$$ where the alternating product of $t$ and $s$ contains $m(s, t) - 1$ factors,
and we have $u(tuw) > tuw$;}
\item{\rm (vii)}{$x$ is weakly complex and there exist $t, u \in S$ with
$m(s, t) = 3$, $ut \ne tu$ and $su = us$ such that $w = sx$ has a reduced
expression of the form $w = utsux'$.}
\endproclaim

\demo{Proof}
Let $\br$ be a reduced expression for $x$ beginning with $s$, and let
$\overline{\br}$ be the set of all reduced expressions for $x$ that are
commutation equivalent to $\br$.

Suppose that some element of $\overline{\br}$ has a reduced expression
beginning with $uu'$, where $u, u'$ are some noncommuting generators in
$S$.  If $u = s$, then we can take $t = u'$ and case (ii) or case (iii)
holds.  If $u \ne s$, then $s$ must commute with both $u$ and $u'$, or it
would not be possible for one element of $\overline{\br}$ to begin with $s$
and another with $uu'$.  This implies that $s$ is distinct from $u$ and $u'$,
and case (iv) applies.

Suppose now that some element of $\overline{\br}$ has a reduced expression
ending with $u'u$, where $u, u'$ are as in the previous paragraph.  
By the arguments in the previous paragraph, 
we may assume that $w$ has a reduced expression
ending in $u'u$, and we are in case (v).

From now on, suppose that neither of the above cases apply.  
This is incompatible with
$x$ being star reducible, so either $x$ is a product of commuting generators,
which is case (i), or $x$ must be weakly complex.  Suppose that the latter
holds.  Now $W$ has Property F, and if $w$ were right star reducible, $x$
would be too.  It must therefore be the case that $w$ has a reduced expression
beginning $uu'$ (where $u$, $u'$ are as before) but that $sw$ has no such
reduced expression.
This means that $s$ must fail to commute with either $u$ or $u'$.  If $s$
fails to commute with $u$, then the earlier analysis shows that case (ii)
or case (iii) applies.  We may now assume that $s$ fails to commute with 
$u'$, and we define $t = u'$.  By Lemma
\secd.5 (i), $x$ has a reduced expression of the form $u w_{st} x'$.
If $m(s, t) > 3$, then $w \in W_c$ has a reduced expression starting
$utstw'$.  Since $tuw \in W_c$ has a reduced expression starting with $st$,
it cannot also have one starting with $u$, so we have $u(tuw) > tuw$; this
is case (vi).  We may now assume that $m(s, t) = 3$, which means that
$w$ has a reduced expression of the form $utsw'$.  If $uw' > w'$ then $w'$,
and hence $sw'$ (because $su = us$) has no reduced expression beginning
with $u$, and case (vi) applies again.  Alternatively, if $uw' < w'$,
then $w$ has a reduced expression of the form $utsux'$, which is case (vii).
\qed\enddemo

\proclaim{Proposition \secf.10}
Suppose $W$ has Property F and Property W, and let $s \in S$ and 
$w \in W_c$.  Writing $$
c_s c_w = \sum_{x \in W_c} \l_x c_x
,$$ we have $\l_x \in \zed^{\geq 0}[\d]$.
\endproclaim

\demo{Proof}
If $sw < w$, this is immediate from Theorem \sece.13, so we may assume that
$sw > w$.  The proof is by induction on $\ell(w)$, the case $\ell(w) = 0$
being trivial.

For the inductive step, we use a case analysis on $x = sw$ based on Lemma 
\secf.9.  In case (i), $x = s_1 s_2 \cdots s_r$ is a product of commuting 
generators, and it is easily verified that $$
c_s c_w = c_x = c_{s_1} c_{s_2} \cdots c_{s_r}
.$$  In cases (ii) and (iii), Proposition \secf.3 (ii) shows that $$
c_s c_w = c_{^*w} + c_{_*w}
,$$ where the star operations are defined with respect to $I = \{s, t\}$, and
as usual, $c_z = 0$ if $z$ is an undefined symbol.

For case (iv), let $I$ be as in the statement of Lemma \secf.9, and write
$w = w_I w^I$.  Let $u$ be as in the statement of Lemma \secf.7.  Then we have $$
c_w = \d^{-1} c_{w_I} c_{uw^I}
.$$  By hypothesis, $s$ commutes with both elements of $I$, which means by
Corollary \secf.5 that $c_s$ commutes with $c_{w_I}$.  We therefore have $$
c_s c_w = \d^{-1} c_{w_I} (c_s c_{uw^I})
.$$  By induction we have $$
c_s c_{uw^I} = \sum_{y \in W_c} \l'_y c_y
,$$ where $\l'_y \in \zed^{\geq 0}[\d]$.  Now $c_u c_{uw^I}
= \d c_{uw^I}$ by Theorem \sece.13, and $c_u$ and $c_s$ commute by hypothesis,
so we must have $$
c_u 
\left( \sum_{y \in W_c} \l'_y c_y \right)
= \d \sum_{y \in W_c} \l'_y c_y
.$$  By Corollary \sece.4, this means that $uy < y$ whenever $\l'_y \ne 0$.
Since $u \in \rdescent{w_I} \cap \ldescent{y}$, we have $$
c_{w_I} c_y = \sum_{z \in W_c} \l''_z c_z
,$$ where $\l''_z \in \d \zed^{\geq 0}$.  Combining these equations completes
the proof in case (iv).

The proof of (v) follows by an argument similar to, but easier than, the 
proof of (iv).

Suppose we are in case (vi), and consider the reduced expression for $w$
given there.  By Proposition \secf.3 (ii), we have $$
c_u c_{uw} = c_w
:$$ the assumption that $u(tuw) > tuw$ implies that ${_*w}$ is undefined with
respect to $I = \{t, u\}$.  Since $s$ commutes with $u$, we have $$
c_s c_w = c_u (c_s c_{uw})
.$$  Although $s(uw) > uw$, we cannot have $suw \in W_c$ because there is
a reduced expression for $suw$ beginning with $w_{st}$.  Using
Proposition \secf.3 (ii) again, we find that $$
c_s c_{uw} = c_{tuw}
,$$ and since $\ell(tuw) < \ell(w)$, we conclude by induction that $$
c_u c_{tuw} = \sum_{x \in W_c} \l_x c_x
,$$ where $\l_x \in \zed^{\geq 0}[\d]$, as required.

Finally, let us suppose that case (vii) holds, and let $I = \{t, u\}$.
Because $t$ fails to commute with both $s$ and $u$, the element
${_*sux'}$ is undefined, and thus (by Proposition \secf.3 (ii) again) we have $$
c_t c_{sux'} = c_{tsux'}
.$$  Now ${^*tsux'} = utsux' = w$ and ${_*tsux'} = sux'$, which implies 
similarly that $$
c_u c_{tsux'} = c_w + c_{sux'}
.$$  This means that $$\eqalign{
c_s c_w &= c_s c_u c_{tsux'} - c_s c_{sux'} \cr
&= c_u (c_s c_{tsux'}) - \d c_{sux'} \cr
&= c_u (c_{sux'}) - \d c_{sux'} \cr
&= \d c_{sux'} - \d c_{sux'} \cr
&= 0,\cr
}$$ where the equalities follow from Theorem \sece.13 and 
Proposition \secf.3 (ii).  This satisfies the hypotheses of the statement 
trivially.
\qed\enddemo

\proclaim{Corollary \secf.11}
Suppose $W$ has Property F and Property W, and let $y, w \in W_c$ be such
that $\ldescent{w} \subsetneq \ldescent{y}$.  Then $\mu(y, w) \geq 0$.
\endproclaim

\demo{Proof}
Let $s \in \ldescent{y} \backslash \ldescent{w}$, so that $sw > w$.  
By Theorem \sece.13, $\mu(y, w)$ is the (integer) coefficient of $c_y$ in 
$c_s c_w$, which by Proposition \secf.10 must be nonnegative.
\qed\enddemo

We return to the issue of positivity of the $\mu(y, w)$ in Corollary \secg.11.

Lemma \secf.8 can now be generalized as follows.

\proclaim{Lemma \secf.12}
Suppose $W$ satisfies Property F and Property W, and let $I = \{s, t\}$ be
a pair of noncommuting generators.  Let $x \in W_c \cap W_I$ and $y \in W_c$.  
Writing $$
c_x c_y = \sum_{y \in W_c} f(x, y, w) c_w
,$$ we have $f(x, y, w) \in \zed^{\geq 0}[\d]$ for all $w$.
\endproclaim

\demo{Proof}
The case $x = 1$ is trivial, so suppose $x \ne 1$ and let $u$ be the unique
element of $\rdescent{x}$.  If $u \in \ldescent{y}$, the claim follows by
Lemma \secf.8, so suppose this is not the case.  Then $$
c_x c_y = \d^{-1} c_x (c_u c_y)
.$$  By Proposition \secf.10, $$
c_u c_y = \sum_{z \in W_c} \l_z c_z
,$$ where $\l_z \in \zed^{\geq 0}[\d]$.  By Theorem \sece.13, $\l_z \ne 0$
implies $uz < z$.  We can now apply Lemma \secf.8 to each term $z$ with
$\l_z \ne 0$ to obtain $$
\d^{-1} c_x c_z = \sum_{x \in W_c} \l'_x c_x
,$$ where each $\l'_x$ lies in $\zed^{\geq 0}[\d]$, and the statement follows.
\qed\enddemo

\proclaim{Theorem \secf.13}
If $W$ satisfies Property F and Property W, then the structure
constants arising from the $c$-basis lie in $\zed^{\geq 0}[\d]$.
\endproclaim

\demo{Proof}
We know that the structure constants lie in $\zed[v, v^{-1}]$, because $TL(X)$
is defined over this ring.  We first note that, as subsets of $\kyu(v)$,
we have $$
\zed^{\geq 0}[\d, \d^{-1}] \cap \zed[v, v^{-1}] = \zed^{\geq 0}[\d]
.$$  Containment in one direction is obvious; to establish the converse,
suppose that $f(v) \in \zed^{\geq 0}[\d, \d^{-1}] \cap \zed[v, v^{-1}]
\backslash \zed^{\geq 0}[\d]$.  Then there is a minimal integer $n > 0$
such that $\d^n f(v) \in \zed^{\geq 0}[\d]$ but $\d^{n-1} f(v) \not\in
\zed^{\geq 0}[\d]$, which means that, as a polynomial in $\d$, $\d^n f(v)$
has a nonzero constant term.
On the other hand, the map $\bar{\ }$ extends to a ring homomorphism of
$\kyu(v)$, and we have $\overline{f(v)} = f(v)$, because $f(v)$ lies in
$\zed^{\geq 0}[\d, \d^{-1}]$.  Since $f(v)$ lies
in the unique factorization domain $\zed[v, v^{-1}]$, $\d^n f(v)$ is an 
$\A$-multiple of the irreducible element $\d$.  Writing $\d^n f(v) = \d g(v)$
and taking images under $\bar{\ }$, we see that $g(v) \in \A$ is 
$\bar{\ }$-invariant.  
However, the $\bar{\ }$-invariant elements of $\A$ are precisely the elements 
of $\zed[\d]$ (because for $k \geq 0$, $\d^k$ is a $\bar{\ }$-invariant
Laurent polynomial with leading term $v^k$)
so in fact $\d^n f(v)$ is a $\zed[\d]$-multiple of $\d$, contradicting
the assumption that $\d^n f(v)$ has nonzero constant term.

It is therefore enough to prove that the structure constants lie in
$\zed^{\geq 0}[\d, \d^{-1}]$.  

Consider a product of two basis elements $c_a c_b$.  We may assume that
$\ell(a), \ell(b) > 1$, or we are done by Lemma \secf.12.
By applying Lemma
\secf.7 repeatedly to each of $c_a$ and $c_b$, we can express $c_a c_b$ as
a finite ordered product of the form $$
\d^{-n} \prod_j c_{w_{I(j)}}
,$$ where for each $j$, $I(j) = \{s_j, t_j\}$ is a pair of noncommuting
generators of $S$ and $\ell(w_{I(j)}) > 0$.  By applying Lemma \secf.12
repeatedly to this product, we find that the structure constants lie in 
$\zed^{\geq 0}[\d, \d^{-1}]$, as required.
\qed\enddemo


It is natural, in the light of the results of \S\sece, to wonder whether the 
$\A$-linear map $\th : \H(W) \ra TL(W)$ satisfying $$
\th(C'_w) = \cases
c_w & \text{ if } w \in W_c, \cr
0 & \text{ otherwise}\cr
\endcases$$ is a homomorphism of algebras.  This is not generally true,
even in the presence of Property F; it fails for example in type $D_4$
\cite{{\bf 15}, Example 2.2.5}.  When the above map is a homomorphism, things 
become much easier, and results such as Theorem \sece.13 are easy to prove.  

The finite Coxeter groups for which $\th$ is a homomorphism were classified
by J. Losonczy and the author in \cite{{\bf 15}}, and for affine Weyl groups by 
Shi in
\cite{{\bf 22}, {\bf 23}}.  The arguments in \cite{{\bf 15}} rely on computer calculations
for types $F_4$, $H_3$ and $H_4$, and the arguments in \cite{{\bf 22}, {\bf 23}}
rely on classification results for Kazhdan--Lusztig cells and on some
deep properties of affine Weyl groups, such as positivity of structure
constants for the $C'$-basis.  It is therefore desirable to find a
conceptual and elementary approach to the problem, which is our aim here.

\proclaim{Proposition \secf.14}
If $W$ has Property S, then $C'_x \in J(X)$ whenever $x \not\in W_c$.
\endproclaim

\demo{Proof}
The proof is by induction on $\ell(x)$, and the base case is vacuous.

If $I \subseteq \ldescent{x}$ or $I \subseteq \rdescent{x}$, then the
argument of the proof of Lemma \secf.2 shows that $C'_x \in J(X)$.

If this is not the case, then $x$ is left or right star reducible to
$x'$, where $x' \not\in W_c$ by Lemma \secd.6.  We treat the case of
left star reducibility, the other case being similar, so write $x
= sx'$.  By Lemma \secf.2, we have $$
C'_s C'_{x'} = C'_x + C'_{_*x} \mod J(X)
.$$  If ${_*x}$ is defined, then ${_*x} \not\in W_c$ by Lemma \secd.6,
and $C'_{_*x} \in J(X)$ by induction.  The same is trivially true if
${_*x}$ is not defined.  Since $C'_{x'} \in J(X)$, the left hand side
of the equation lies in $J(X)$.  It follows that $C'_x \in J(X)$, as
required.
\qed\enddemo


\proclaim{Corollary \secf.15}
If $W$ has Property S, then $\te_w \in v^{-1}\L$ for all complex $w \in 
W$.  In particular, $W$ has Property W.
\endproclaim

\demo{Proof}
This follows from Proposition 6.14 and the equivalence of parts (ii) and (v)
of \cite{{\bf 15}, Theorem 2.2.3.}
\qed\enddemo


\proclaim{Theorem \secf.16}
Suppose that $W$ has Property F and Property S.  Let $x, y \in W$ and
write $$
C'_x C'_y = \sum_{z \in W} g(x, y, z) C'_z
.$$  If $z \in W_c$, then $g(x, y, z) \in \zed^{\geq 0}[\d] \subset
\zed^{\geq 0}[v, v^{-1}]$.
\endproclaim

\demo{Proof}
Applying $\th$ to the equation in the statement and using Proposition
\secf.3 (i) and Proposition \secf.14, we obtain $g(x, y, z) = 0$ unless 
$x, y \in W_c$, and in the latter case, we have $$
c_x c_y = \sum_{z \in W_c} g(x, y, z) c_z
.$$  The result now follows from Theorem \secf.13.
\qed\enddemo

\head \secg.  Computing the $\mu(x, w)$ using generalized Jones traces 
\endhead


The main aim of \S\secg\  is to show how, in many cases, the coefficients 
$\mu(y, w)$, for $y, w \in W_c$, may be computed nonrecursively using 
a(ny) generalized Jones trace.  To the best of our knowledge, this result
is new even in type $A$.

To this end, we need some combinatorial lemmas involving fully commutative
elements.

\definition{Definition \secg.1}
Let $W$ be any Coxeter group and let $w \in W_c$.  We define $n(w)$ to
be the maximum integer $k$ such that $w$ has a reduced expression of the
form $w = w_1 w_2 w_3$, where $\ell(w_2) = k$ and $w_2$ is a product of
commuting generators.
\enddefinition

The following result was proved by Shi \cite{{\bf 24}, Lemma 2.9} for finite and
affine Weyl groups, but it is an easy exercise to prove it for arbitrary
Coxeter groups.

\proclaim{Lemma \secg.2}
Let $W$ be any Coxeter group and let $w \in W_c$.  If $w$ is left
(or right) star reducible to $x \in W_c$, then $n(x) = n(w)$.
\qed\endproclaim

By iterating Lemma \secg.2, we obtain the following 

\proclaim{Corollary \secg.3}
Suppose $W$ has Property F, and let $w \in W_c$.  Then $w$ is star
reducible to a product of $n(w)$ generators.
\qed\endproclaim

\proclaim{Lemma \secg.4}
Suppose $w \in W_c$ is such that
$|\ldescent{w}| = n(w)$ (respectively, $|\rdescent{w}| = n(w)$).  Then
if $w$ is left (respectively, right) star reducible to $x$, we
have $|\ldescent{w}| = |\ldescent{x}|$ and $\rdescent{w} = \rdescent{x}$
(respectively, $|\rdescent{w}| = |\rdescent{x}|$ and 
$\ldescent{w} = \ldescent{x}$).
\endproclaim

\demo{Proof}
We deal with the case where $|\ldescent{w}| = n(w)$, the other case being
similar.  It is immediate from the definitions that if $y \in W_c$ is left
star reducible to $y'$, then $|\ldescent{y'}| \geq |\ldescent{y}|$ and
$\rdescent{y'} = \rdescent{y}$.  The definition of $n(y)$ shows that we
always have $\max\{|\ldescent{y}|, |\rdescent{y}|\} \leq n(y)$.  Lemma
\secg.2 and the hypothesis $|\ldescent{w}| = n(w)$ thus force equality 
as required.
\qed\enddemo

\definition{Definition \secg.5}
Suppose that the Coxeter graph $X$ is bipartite, and let $$\e : S \ra 
\{0, 1\}$$ be a labelling of $S$ corresponding to a 2-colouring of the 
graph.  If $J \subset S$ is a subset of commuting
generators, we define $$
k_\e(J) = (-1)^{|J \cap \e^{-1}(0)|}
.$$  For $w \in W_c$, we define $k_\e(w) \in \{\pm 
1\}$ by $$
k_\e(w) = k_\e(\ldescent{w}) \times k_\e(\rdescent{w})
.$$
\enddefinition

\proclaim{Lemma \secg.6}
Let $W$ be a Coxeter group with $X$ bipartite and $\e$ as in
Definition \secg.5.  
Let $w \in W_c$ be such that $|\ldescent{w}| = n(w)$ (respectively, 
$|\rdescent{w}| = n(w)$), and suppose $w$ is left (respectively, right)
star reducible to $x \in W_c$.  Then $k_\e(w) = -k_\e(x)$.
\endproclaim

\demo{Proof}
By symmetry, we only deal with the case of left star reducibility.  If
$w$ is left star reducible to $x$ with respect to $I = \{s, t\}$, then
$\e(I) = \{0, 1\}$.  It follows from Lemma \secg.4 that 
$k_\e(\ldescent{w}) = -k_\e(\ldescent{x})$ and $k_\e(\rdescent{w}) = 
k_\e(\rdescent{x})$, and the claim follows.
\qed\enddemo

\proclaim{Lemma \secg.7}
Let $W$ be a Coxeter group with $X$ bipartite and $\e$ as in
Definition \secg.5, and suppose also that $W$ has Property F.  Let
$w \in W_c$ be such that $\ldescent{w} = \rdescent{w}$ is a set of size
$n(w)$.  Then $\ell(w) = n(w) \mod 2$.
\endproclaim

\demo{Proof}
Choose a function $\e$ as in Definition \secg.5.
The hypothesis that $\ldescent{w} = \rdescent{w}$ means that $k_\e(w) = 1$.
By Corollary \secg.3, $w$ is star reducible to a product $y$ of $n(w)$ 
generators; since $\ldescent{y} = \rdescent{y}$, we have $k_\e(y) = 1$
as well.  By Lemma \secg.6, there must have been an even number of star
operations applied to reduce $w$ to $y$, each of which decreases the length
by $1$.  The claim now follows.
\qed\enddemo

We now turn our attention to Coxeter groups having Property B.  It is clear
from Hypothesis \secb.2 (ii) that if $x, y \in W_c$ are distinct, then
$v \lrt{\te_x}{\te_y} \in \A^-$.  We will show that in many important cases,
we in fact have $v \lrt{\te_x}{\te_y} \in v^{-1} \A^-$.

\proclaim{Lemma \secg.8}
Suppose that $W$ has Property B, let $x, y \in W_c$ be distinct
elements, and let $f(v) = v \lrt{\te_x}{\te_y}$.  
\item{\rm (i)}
{If $\e_x = \e_y$, then $f(v) \in v^{-1} \A^-$.}
\item{\rm (ii)}
{If $x^{-1} y \in S$ or $y x^{-1} \in S$, then $f(v) \in v^{-1} \A^-$.}
\item{\rm (iii)}
{If $\ldescent{x} \ne \ldescent{y}$ or $\rdescent{x} \ne \rdescent{y}$, then
$f(v) \in v^{-1} \A^-$.}
\endproclaim

\demo{Proof}
We assume, by Lemma \secb.8, that the form $\lrt{\ }{\ }$ is homogeneous.
This means that $\lrt{\te_x}{\te_y} \in \zed[v^{-2}]$ if $\e_x = \e_y$,
and  $\lrt{\te_x}{\te_y} \in v^{-1} \zed[v^{-2}]$ otherwise.  If we are
in the former case and $x \ne y$, we have $\lrt{\te_x}{\te_y} \in 
v^{-2} \zed[v^{-2}]$, and (i) follows.

To prove (ii), let us assume that $x = uy < y$ for some $u \in S$; the other
case is similar.  By Hypothesis \secb.2 (ii), we have $$\eqalign{
1 
&= \lrt{\te_{ux}}{\te_{ux}} \mod v^{-1} \A^- \cr
&= \lrt{\te_x}{\te_u \te_{ux}} \mod v^{-1} \A^- \cr
&= \lrt{\te_x}{\te_{x}} + 
(v - v^{-1}) \lrt{\te_x}{\te_{ux}} \mod v^{-1} \A^- \cr
&= 1 + v \lrt{\te_x}{\te_{ux}} \mod v^{-1} \A^-, \cr
}$$ which shows that $v \lrt{\te_x}{\te_y} \in v^{-1} \A^-$, as required.
(Note that, for the second equality, we have $\te_{ux} = \te_u \te_x$ because
$ux > x$.)

For (iii), let us assume that $\ldescent{y} \not\subseteq \ldescent{x}$;
the other cases follow similarly.  (Recall that $\lrt{\te_x}{\te_y}
= \lrt{\te_y}{\te_x}$.)  Let $u \in \ldescent{y} \backslash \ldescent{x}$.
We may assume that $x \ne uy$ or we are
done by part (ii).  Using the identity $$
v \te_y = \te_u \te_y + v^{-1} \te_y - \te_{uy}
,$$ we have $$\eqalign{
v \lrt{\te_x}{\te_y} 
&= \lrt{\te_x}{\te_u \te_y} + v^{-1} \lrt{\te_x}{\te_y} - 
\lrt{\te_x}{\te_{uy}} \cr
&= \lrt{\te_x}{\te_u \te_y} \mod v^{-1} \A^- \cr
&= \lrt{\te_{ux}}{\te_y} \mod v^{-1} \A^- \cr
&= 0 \mod v^{-1} \A^-,\cr
}$$ as required.
\qed\enddemo

\proclaim{Proposition \secg.9}
Let $W$ be a Coxeter group with Property B and Property F such that
the graph $X$ is bipartite.  If the bilinear form is homogeneous, 
then for $x, y \in W_c$ we have $$
\lrt{\te_x}{\te_y} = 
\cases 1 \mod v^{-2} \A^- & \text{ if } x = y,\cr
0 \mod v^{-2} \A^- & \text{ otherwise.} \cr
\endcases$$  In other words, for any distinct elements $x, y \in W_c$,
we have $$
v \lrt{\te_x}{\te_y} \in v^{-1} \A^-
.$$
\endproclaim

\demo{Proof}
The proof is by induction on $n = \ell(x) + \ell(y)$.  By Lemma \secg.8 (i),
we only need deal with the case where $n$ is odd.  The base case is then
$n = 1$, which says that $$
v \lrt{\te_s}{\te_1} \in v^{-1} \A^-
,$$ where $s \in S$.  This also follows from Lemma \secg.8 (ii).

Suppose now that $n = k$ for some odd number $k$, and that
the statement is known to be true for all $n < k$.  By Lemma \secg.8 (iii),
we may assume that $\ldescent{x} = \ldescent{y}$ and 
$\rdescent{x} = \rdescent{y}$.

Suppose at first that $x$ is not the product of commuting generators.  By
Property F, $x$ is either left or right star reducible; we only treat the
case of left star reducibility by symmetry.  In this case, there exist
noncommuting generators $s, t$ such that $x = stx'$ and $y = sy'$ are
reduced.  By Lemma \secc.7 (i) and the inductive hypothesis, we have 
$$\eqalign{
f(v) &= v \lrt{\te_s \te_t \te_{x'}}{\te_s \te_{y'}} \cr
&= v \lrt{\te_{tx'}}{\te_t \te_{y}} 
+ v \lrt{\te_{x'}}{\te_y} 
- v \lrt{\te_{tx'}}{\te_{y'}} \cr
&= v \lrt{\te_{tx'}}{\te_t \te_{y}} \mod v^{-1} \A^-.
}$$  Since $sy < y$ and $y \in W_c$, we must have $ty > y$.  If $ty \not\in
W_c$, Lemma \secd.5 (i) shows that $ty$ has a reduced expression beginning in
$w_{st}$.  In this case, Lemma \secd.11 (ii) shows that $$
v \lrt{\te_{tx'}}{\te_t \te_{y}} =
- \lrt{\te_{tx'}}{\te_{y} + c \te_{sty}}
\mod v^{-1} \A^-
,$$ where $c = 1$ if $sty \in W_c$, and $c = 0$ otherwise.  If the above
expression does not lie in $v^{-1} \A^-$, we must have either $y = tx'$ or
both $sty \in W_c$ and $sty = tx'$.  The former situation is impossible
because $sy < y$ and $stx' > tx'$.  The latter situation also cannot occur,
because it implies that $x = stx' = ty$, which contradicts $x \in W_c$ and
$ty \not\in W_c$.  We conclude that in fact $ty \in W_c$.  In summary,
what we have shown is that, with respect to $I = \{s, t\}$, we have 
$$\eqalignno{
v \lrt{\te_x}{\te_y} &= v \lrt{\te_{_*x}}{\te_{^*y}} \mod v^{-1} \A^-
, & (11)}$$ where we interpret $\te_{^*y}$ as $0$ if ${^*y}$ is not defined.

We can now apply (11) (and its right-handed version) repeatedly, which
will either prove the claim along the way or result in consideration of
a quantity $v \lrt{\te_x}{\te_y}$, where $\ell(x) + \ell(y) = k$,
$x$ is a product of $a$ commuting generators and $\ldescent{y} = \rdescent{y}$
consists of the same $a$ commuting generators.  From the definition of
$n(y)$, we see that $n(y) \geq a$.  If we have $n(y) > a$, we can exchange
the roles of $x$ and $y$ and again apply (11) (and its right-handed version) 
repeatedly until this is no longer possible.  If this does not prove the claim
along the way, Corollary \secg.3 shows that we obtain a quantity 
$v \lrt{\te_{x'}}{\te_{y'}}$, where $\ell(x') + \ell(y') = k$, $x'$ is a 
product of $n(y) > a$ commuting
generators and $\ldescent{y'} = \rdescent{y'}$
consists of the same $n(y)$ commuting generators.  If we still have $n(y')
> n(y)$, we can repeat the same process; eventually this must terminate because
the $n$-values strictly increase at each step, and they are bounded above
by $k$.  

We have now reduced consideration to the case of $v \lrt{\te_x}{\te_y}$,
where $x$ is a product of $n(y)$ commuting generators, and $\ldescent{y}
= \rdescent{y}$ consists of the same $n(y)$ commuting generators.  Since
$X$ is bipartite, Lemma \secg.7 now applies to show that $\e_x = \e_y$, and
the proof is completed by Lemma \secg.8 (i).
\qed\enddemo

We may now prove the main result of this section.

\proclaim{Theorem \secg.10}
Let $W$ be a Coxeter group with Property B and Property F such that
the graph $X$ is bipartite, and assume that the form $\lrt{\ }{\ }$ is
homogeneous.  Then for any elements $x, y \in W_c$,
the coefficient of $v^{-1}$ in $\lrt{c_x}{c_y}$ is $\tmu(x, y)$.
\endproclaim

\demo{Proof}
Without loss of generality, we suppose that $\ell(y) \geq \ell(x)$.
By equation (2), we have $$
v \lrt{c_x}{c_y} = 
v \lrt
{\sum_{a \in W_c} p^*(a, x) \te_a}
{\sum_{b \in W_c} p^*(b, y) \te_b}
.$$  Recall that $p^*(c, d) \in v^{-1} \A^-$ unless $c = d$, and by Theorem
\sece.13, the coefficient of $v^{-1}$ in $p^*(c, d)$ is $\mu(c, d)$.
Proposition \secg.9 shows that $v \lrt{\te_a}{\te_b} \in v^{-1} \A^-$ unless 
$a = b$.

It follows that the only way we can have $$
v \lrt{p^*(a, x) \te_a}{p^*(b, y) \te_b} \not\in v^{-1} \A^-
$$ is if both $a = b$, and either $a = x$ or $b = y$ (or both).  However,
if $a = b$ and $a = x$ and $b = y$, then $\e_x = \e_y$ and $\tmu(x, y) = 0$, 
and the coefficient of $v^{-1}$ in $\lrt{c_x}{c_y}$ is zero by homogeneity, 
which completes the proof.  If $a = b$ and $b = y$ but $a \ne x$, we may
assume that $a < x$, which means that $\ell(a) < \ell(x) \leq \ell(y)$,
contradicting $a = b$.  The only case left to consider is when $a = b$,
$a = x$ and $b \ne y$.  In this case, we have $$
v \lrt{p^*(a, x) \te_a}{p^*(b, y) \te_b} 
= v \lrt{\te_x}{p^*(x, y) \te_x} = \mu(x, y) \mod v^{-1} \A^-
,$$ as required.
\qed\enddemo

\proclaim{Corollary \secg.11}
If $W$ is a Coxeter group with Property B and Property F such that
the graph $X$ is bipartite, and such that the trace $\t$ is homogeneous
and positive (in the sense of Definition \secb.9), then the integers 
$\tmu(x, y)$ are nonnegative.
\endproclaim

\demo{Proof}
By Theorem \secf.13, the product
$c_x c_{y^{-1}}$ is a $\zed^{\geq 0}[\d]$-linear combination of basis
elements.  Since $\t$ is positive, we have $$
\lrt{c_x}{c_y} = \t(c_x c_{y^{-1}}) \in \zed^{\geq 0}[v, v^{-1}]
,$$ and the result follows from Theorem \secg.10.
\qed\enddemo

\remark{Remark \secg.12}
Note that in the simply laced case, Corollary \secg.11 is obvious from 
Proposition \sece.9, which has at most one nonzero term on each side of the 
equation.  (In fact, in this case, it is clear that the $\tmu(x, y)$ are
all equal to $0$ or $1$.)  In the case of type $ADE$, 
Graham \cite{{\bf 7}, proof of Theorem 9.9} gives a nice characterization of 
those $x \in W_c$ for which $x \leq w$ for some fixed $w \in W_c$: such $x$
arise from the basis elements $c_x$ obtained by deleting a single generator 
from the monomial $c_w$.
It is not clear if this could be generalized to non-simply-laced cases.
However, given elements $x, w \in W_c$, Graham's method for computing
$\mu(x, w)$ is recursive, unlike Theorem \secg.10 above.
\endremark

\remark{Remark \secg.13}
Closed formulae for Kazhdan--Lusztig polynomials have been developed by
Brenti \cite{{\bf 2}}; these involve taking the sum over certain chains.
However, when an explicit construction for the trace $\t$ is known, 
Theorem \secg.10 typically requires very little computation indeed, as we
illustrate below.  This means that one can be very explicit about the
values $\tmu(x, w)$; for example, one can show using diagram calculus
methods in \cite{{\bf 10}} that in type $B$ or type $H_n$ (even when $n$ is
arbitrarily large), the integers $\tmu(x, w)$ are always $0$ or $1$ when
$x, w \in W_c$.  It would be interesting to know if this holds generally.
\endremark

\remark{Remark \secg.14}
The hypothesis that $X$ be bipartite cannot be removed from Theorem \secg.10.
For example, in type $\widehat{A}_2$, which does satisfy Property B and 
Property F, it is possible to find a homogeneous bilinear form 
$\lrt{\ }{\ }$ such that $$
\lrt{c_x}{c_y} = N
$$ where $S = \{s_1, s_2, s_3\}$, $x = s_1$, $y = s_1 s_2 s_3 s_1$,
and any given integer $N$.

It is possible to prove Theorem \secg.10 for some Coxeter groups that do not
have Property F, such as type $\widehat{A}_n$ for $n$ odd, but this requires
significant modifications to the arguments.
\endremark

To the best of our knowledge, Theorem \secg.10 is new even in type $A$.  In this
case, the result shows how Jones' trace on the Temperley--Lieb algebra may
be used to compute all values $\mu(x, w)$ for which $x, w \in W_c$.

\example{Example \secg.15}
Let $W$ be a Coxeter group of type $A_3$, and let $\t$ be the homogeneous
trace of Remark \secb.5 (i).  Let $x = s_2$ and $y = s_2 s_1 s_3 s_2$, where
the generating set $S$ is indexed in the obvious way.  Using the 
Temperley--Lieb diagram calculus, we see immediately from Figure 1 
that the diagram corresponding to $$
\t(c_x c_{y^{-1}}) = \t(c_{s_2} c_{s_2} c_{s_3} c_{s_1} c_{s_2})
$$ has $3$ closed loops, and so we have $$
\t(c_x c_{y^{-1}}) = v^{-4} (v + v^{-1})^3
,$$ in which the coefficient of $v^{-1}$ is $1$.  This proves that
$\mu(x, y) = 1$.  Since $P_{x, y}(q)$ has degree at most $1$ and constant
term $1$, this recovers the well-known result that $P_{x, y} (q) = 1 + q$
for these elements.
\endexample

\topcaption{Figure 1} Computation of $\mu(s_2, s_2 s_1 s_3 s_2)$
\endcaption
\centerline{
\hbox to 1.361in{
\vbox to 1.930in{\vfill
        \includegraphics{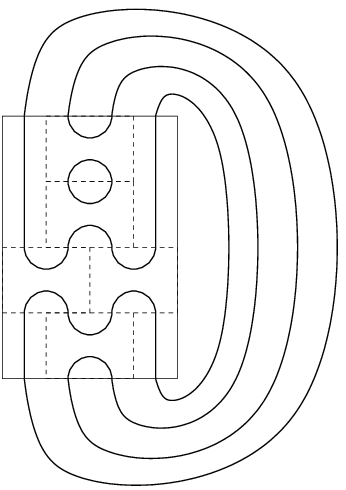}
}
\hfill}
}

\head \secconc.  Overview and conclusion \endhead

In the sequel \cite{{\bf 11}} to this paper, we show how Property W is 
in fact a consequence of Property F \cite{{\bf 11}, Theorem 4.6 (i)}.  
As a by-product, we show in \cite{{\bf 11}, Theorem 4.3}
how, under this hypothesis,
we have $\te_w \in \L$ for all complex $w \in W$, or, equivalently (if
Property B holds), $$
\lrt{\te_x}{\te_w} \in \A^- \text{ for all } x, w \in W
.$$  This result is one of the 
``projection properties'' studied in \cite{{\bf 14}, {\bf 19}}.
It is obvious if Property S holds, but is nontrivial otherwise, for example
in the case of type $D$, where it was proved by Losonczy
\cite{{\bf 19}}.

A main theme of the papers \cite{{\bf 15}, {\bf 22}, {\bf 23}} is the compatibility between
Kazh\-dan--Lusz\-tig cells and fully commutative elements.  
In terms of Property B, this asks whether $$
\lrt{\te_x}{\te_w} \in \cases
\A^- \text{ for all } x, w \in W \text{ and} \cr
v^{-1} \A^- \text{ if } x \not\in W_c \text{ or } w \not\in W_c. \cr
\endcases
$$  The results of this paper allow more elegant proofs of these results.
In particular, \cite{{\bf 23}, Lemma 2.4}, which relies on the theory of
cells in affine Weyl groups, becomes unnecessary due to Proposition 6.14.
It is also possible to apply Property S to avoid the ad hoc arguments in
\cite{{\bf 23}, Appendix} based on cell classifications.

It would be interesting to know whether generalized Jones traces exist 
for all Coxeter systems, but it seems likely that an elementary proof of 
this would be difficult.

\head Acknowledgements \endhead

I thank J. Losonczy for many helpful comments on an early version of this
paper.  I am also grateful to the referee for reading the paper carefully
and suggesting some improvements.

\leftheadtext{} \rightheadtext{}

\end